\documentclass[journal]{IEEEtran}
\IEEEoverridecommandlockouts
\usepackage{amsmath,amsfonts,amssymb,euscript, graphicx, 
epsfig,
enumerate,float,afterpage, subfigure, ifthen}%

\newtheorem{thm}{Theorem}

\newtheorem{lem}{Lemma}
\newtheorem{defn}{Definition}

\newcommand{\expect}[1]{\mathbb{E}\left\{#1\right\}}
\newcommand{\defequiv}{\mbox{\raisebox{-.3ex}{$\overset{\vartriangle}{=}$}}}

\newcommand{\bv}[1]{{\boldsymbol{#1} }}
\newcommand{\script}[1]{{{\cal{#1} }}}

\begin{document}

\title
  {Dynamic Optimization and Learning for Renewal Systems}
\author{Michael J. Neely\\$\vspace{-.4in}$%University of Southern California% \\ http://www-rcf.usc.edu/$\sim$mjneely%
\thanks{Michael J. Neely is with the  Electrical Engineering department at the University
of Southern California, Los Angeles, CA.}
\thanks{This paper was presented in part at the Asilomar Conference on Signals, Systems, and Computers, 
Pacific Grove, CA, 2010 \cite{renewals-asilomar2010}. 
This material is supported in part  by one or more of 
the following: the DARPA IT-MANET program
grant W911NF-07-0028, 
the NSF Career grant CCF-0747525, and continuing through participation in the 
Network Science Collaborative Technology Alliance sponsored
by the U.S. Army Research Laboratory.}}

\markboth{Extended Version}{Neely}

\maketitle

%UNCOMMENTING THIS REMOVES PAGE NUMBERS
%\thispagestyle{empty}
%\pagestyle{empty} 

\begin{abstract}   
We consider the problem of optimizing time averages
in systems with independent and identically distributed 
behavior over renewal frames.
This includes scheduling and task processing to maximize utility in stochastic
networks  with 
variable length scheduling modes.
Every frame, a new policy is implemented that affects the frame
size and that creates a vector of attributes.Ê An algorithm
is developed for choosing policies on each frame in order to
maximize a concave function of the time average attribute
vector, subject to additional time average constraints.
The algorithm is based on Lyapunov optimization concepts
and involves minimizing a ``drift-plus-penalty'' ratio over each frame. 
The algorithm can learn
efficient behavior without a-priori statistical knowledge by
sampling from the past.  Our framework is applicable to a large
class of problems, including Markov decision problems.
\end{abstract}

\section{Introduction} 

Consider a stochastic system that regularly experiences times when the system
state is refreshed, called \emph{renewal times}.  The goal is to develop a control algorithm 
that maximizes the time
average of a reward process associated with the system, subject to time average constraints
on a collection of penalty processes.  
The renewal-reward theorem is a simple and elegant technique for computing time averages
in such systems (see, for example, \cite{gallager}\cite{ross-prob}).  However, the renewal-reward theorem requires 
random events to be independent and identically distributed (i.i.d.) over each renewal frame.  While this i.i.d. assumption
may hold if a single control law is implemented repeatedly, it is often difficult to choose in advance a
single control law that optimizes the system subject to the desired constraints. 
This paper investigates the situation where the control policies used may differ from frame to frame, and are
designed to dynamically solve the problem of interest. 

This renewal problem arises in many different applications.  One application of interest is a 
\emph{task processing 
network}.  For example, consider 
a network of wireless devices that repeatedly collaborate to accomplish tasks (such
as reporting sensor data to a destination, or performing distributed computation on data). Tasks are performed
one after the other, and for each task we must decide what modes of operation and communication to use, possibly allowing 
some nodes of the network to remain idle to save power.  It is then important to make decisions that 
maximize the time average utility associated with task processing, subject to time average power constraints at each node. 
Alternatively, one may want to minimize time average power, subject to constraints on utility and on the ``left-over'' 
communication rates available for data that is not associated with the task processing. 

This paper develops a general framework for solving
such problems.  To do so, we extend the theory of Lyapunov optimization from \cite{now}. 
Specifically, work in \cite{now} considers
discrete time queueing networks and develops 
a simple \emph{drift-plus-penalty} rule for making optimal decisions.  These
decisions are made in a greedy manner every slot based only on 
the observed traffic and channel conditions for that slot, without requiring 
a-priori knowledge of the underlying probability distribution. 
However, the work in \cite{now} assumes
all slots have fixed length, the random network condition is observed at the beginning
of each slot and does not change over the slot,  and this condition
is not  influenced by control actions. 
The general renewal problem treated in the current paper
is more complex because each frame may have a different length and
may contain a sequence of random events.  The frame length and the random event sequence may depend on the 
control  decisions  made over the course of the frame. Rather than making a single decision every slot, every frame we must
specify a \emph{policy}, being a contingency plan for making decisions over the course of the frame in reaction to the 
resulting system events. 

This paper solves the general problem with a conceptually simple technique that chooses a policy to minimize
a  \emph{drift-plus-penalty ratio} every
frame. 
We first develop algorithms for minimizing the time average of a penalty process subject to a collection of time average
constraints.  We then consider maximization of a concave function of a vector of time average attributes
subject to similar constraints.   This utility maximization problem is challenging because of the variable frame length.  We overcome
this challenge with a novel transformation together with a variation of Jensen's inequality. 

While this paper focuses on task processing applications, we note that our renewal framework can also handle
\emph{Markov decision problems}.  Specifically, suppose the system operates according to either a continuous or discrete time
Markov chain with control-dependent transition probabilities.  If the chain has a recurrent state, then renewals can be defined
as re-visitations to this state, and the same drift-plus-penalty ratio technique can be applied.  However,  the
drift-plus-penalty ratio may be difficult to optimize for Markov decision problems with high dimension
(see also  \cite{neely-mdp-cdc09}). 

Prior work on learning algorithms for Markov decision problems is in 
\cite{self-learning-mdp}, and related work in \cite{q-learning-mimo}\cite{energy-delay-approx}\cite{mihaela-dp1}\cite{mihaela-dp2}
considers learning for optimization of energy and 
delay in queueing systems.  The works \cite{self-learning-mdp}-\cite{mihaela-dp2}
use stochastic approximation 
theory and two-timescale convergence analysis.  The Lagrange multiplier updates 
in \cite{self-learning-mdp}-\cite{mihaela-dp2} are analogous
to the \emph{virtual queue updates} we use in this paper.  However, the Lyapunov optimization framework we use is different and
does not require a two-timescale approach.  It also 
provides more explicit bounds on convergence times and 
deviations from optimality, and allows a broader class of problems such as
task processing problems. 

The Lyapunov optimization technique that we use in this paper is based on our previous work in 
\cite{now}\cite{neely-thesis}\cite{neely-energy-it}\cite{neely-fairness-infocom05} that develops the drift-plus-penalty
method for stochastic network optimization, including opportunistic scheduling for throughput-utility maximization \cite{now}\cite{neely-thesis}\cite{neely-fairness-infocom05} and average power 
minimization \cite{neely-energy-it} (see also \cite{sno-text}).  Alternative ``fluid-based'' 
stochastic optimization techniques for queueing networks are developed in 
\cite{atilla-fairness-ton}\cite{stolyar-greedy}\cite{stolyar-gpd-gen}\cite{primal-dual-cmu}, and dual and primal-dual algorithms
for systems without queues, based on tracking a corresponding static optimization problem, 
are in \cite{lin-shroff-cdc04}\cite{vijay-allerton02}\cite{prop-fair-down}.
Our current paper considers the more complex renewal problem, and 
leverages ideas in \cite{neely-mdp-cdc09}\cite{chihping-utility-round-robin}, where
\cite{neely-mdp-cdc09} considers a frame-based Lyapunov framework for Markov decision problems involving
network delay, and \cite{chihping-utility-round-robin} develops a ratio rule for utility optimization in wireless systems
with variable length
frames and time-correlated channels. 

%An important class of task processing problems are recently treated in  \cite{network-corroboration-arxiv}.  That work
%considers a system 
%where multiple wireless ``reporting nodes'' select data formats (e.g., ``voice'' or ``video'')  
%in which to deliver sensed information.  The work 
%in \cite{network-corroboration-arxiv} also uses a renewal structure.  However, it assumes 
%a single random event occurs at the beginning of each renewal frame, and the event and frame size are not influenced by control actions.  
%More general problems can be treated using the theory developed in the current
%paper. 

Recent work in \cite{network-corroboration-arxiv} considers a task processing system where multiple
wireless ``reporting nodes'' select data formats (e.g., ``voice'' or ``video'')  
in which to deliver sensed information.  The work 
\cite{network-corroboration-arxiv} also uses a renewal structure.  However, it assumes 
a single random event occurs at the beginning of each renewal frame, and the event and frame size are not influenced by control actions.  
More general problems can be treated using the theory developed in the current
paper.

\section{Renewal System Model}

%\begin{figure}[htbp]
\begin{figure}[t]
   \centering
   \includegraphics[width=3.5in, height=0.75in]{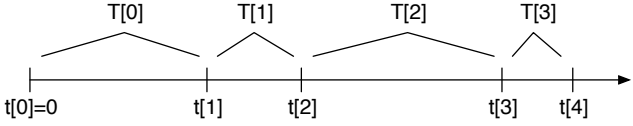} % requires the graphicx package
   \caption{A timeline illustrating renewal frames for the system.}
   \label{fig:renewal-timeline}
\end{figure}

Consider a system that operates over \emph{renewal frames}.  Specifically, consider the timeline of non-negative real 
times $t \geq 0$, and suppose this timeline is segmented into successive frames of duration $\{T[0], T[1], T[2], \ldots\}$, as shown
in Fig. \ref{fig:renewal-timeline}.  Define $t[0]=0$, and for each positive integer $r$ define $t[r]$ as the \emph{$r$th renewal time}: 
\begin{eqnarray*}
&t[r] \defequiv \sum_{i=0}^{r-1} T[i]&
\end{eqnarray*}
The interval of all times $t$ such that $t[r] \leq t < t[r+1]$ is defined as the \emph{$r$th renewal frame}, 
defined for each $r \in \{0, 1, 2, \ldots\}$. 

At the beginning of each renewal frame $r$, the controller selects a \emph{policy} $\pi[r]$ from an abstract policy space
$\script{P}$, and implements the policy over the duration of the frame.   There may be random events that arise over the renewal 
frame (with distributions that are possibly dependent on the policy), and the policy specifies a contingency plan for reacting
to these events. 
The policy incurs a vector of \emph{penalties} $\bv{y}[r] = (y_0[r], y_1[r], \ldots, y_L[r])$ and \emph{attributes} $\bv{x}[r] = (x_1[r], \ldots, x_M[r])$ for some integers $L\geq0$, $M\geq0$ (where $L=0$ corresponds to problems without $\bv{y}[r]$ penalties, 
and $M=0$ corresponds to problems without $\bv{x}[r]$ attributes). 
  The policy may
also affect the renewal frame duration $T[r]$.  Formally, the values $T[r]$, $y_l[r]$, $x_m[r]$ are determined by 
\emph{random functions} $\hat{T}(\cdot)$, $\hat{y}_l(\cdot)$, $\hat{x}_m(\cdot)$ of the policy $\pi[r]$: 
\begin{eqnarray}
T[r] &\defequiv& \hat{T}(\pi[r]) \label{eq:functions1} \\
y_l[r] &\defequiv& \hat{y}_l(\pi[r]) \: \: \forall l \in \{0, 1, \ldots, L\} \label{eq:functions2} \\
x_m[r] &\defequiv& \hat{x}_m(\pi[r]) \: \: \forall m \in \{1, \ldots, M\} \label{eq:functions3} 
\end{eqnarray}

We assume  the values of $[\hat{T}(\pi[r]), (\hat{y}_l(\pi[r])), (\hat{x}_m(\pi[r]))]$ for frame $r$ are conditionally 
independent of events in previous frames given the particular policy $\pi = \pi[r]$, and are identically distributed over all frames that use
the same policy $\pi$.  

Consider now a particular control algorithm that chooses policies $\pi[r] \in \script{P}$ every frame $r$ according to some well
defined (possibly probabilistic) rule, and define the following
frame-average expectations, defined for integers $R>0$: 
\begin{eqnarray}
\overline{T}[R] \defequiv \frac{1}{R}\sum_{r=0}^{R-1} \expect{T[r]} \: \: , \: \: 
\overline{y}_l[R] \defequiv \frac{1}{R}\sum_{r=0}^{R-1}\expect{y_l[r]} \label{eq:yav-r-frame} 
\end{eqnarray}
where we recall that $T[r]$, $y_l[r]$, $x_m[r]$ depend on the policy $\pi[r]$ by (\ref{eq:functions1})-(\ref{eq:functions3}). 
Define $\overline{x}_m[R]$ similarly, and 
define the infinite horizon frame-average expectations $\overline{T}$, $\overline{y}_l$, $\overline{x}_m$  by: 
\[ (\overline{T}, \overline{y}_l, \overline{x}_m)  = \lim_{R\rightarrow\infty} (\overline{T}[R], \overline{y}_l[R], \overline{x}_m[R])\]
where we temporarily assume the limits are well defined. 

\subsection{Optimization Objective} 
The first type of problem we consider uses only penalties $\bv{y}[r]$:  We must choose a policy $\pi[r] \in \script{P}$ every frame
$r$  to minimize the ratio $\overline{y}_0/\overline{T}$ subject to constraints on $\overline{y}_l/\overline{T}$: 
\begin{eqnarray}
\mbox{Minimize:} & \overline{y}_0/\overline{T} \label{eq:p1} \\
\mbox{Subject to:} & \overline{y}_l/\overline{T} \leq c_l \: \: \forall l \in \{1, \ldots, L\} \label{eq:p2}  \\
& \pi[r] \in \script{P} \: \: \forall r \in \{0, 1, 2, \ldots\} \label{eq:p3} 
\end{eqnarray}
where $c_l$ for $l \in \{1, \ldots, L\}$ are a given collection of real-valued (possibly negative) constants. 

The motivation for looking at the ratio $\overline{y}_l/\overline{T}$ is that it defines the \emph{time average penalty
associated with the $y_l[r]$ process}.  To see this, suppose the following limits converge to constants $y_l^{av}$
and $T^{av}$ with probability 1: 
\begin{eqnarray*}
 \lim_{R\rightarrow\infty} \frac{1}{R}\sum_{r=0}^{R-1} y_l[r] = y_l^{av} \: \: , \: \: 
 \lim_{R\rightarrow\infty} \frac{1}{R}\sum_{r=0}^{R-1} T[r] = T^{av} \: \: (w.p.1) 
\end{eqnarray*} 
Under very mild conditions, 
the existence of the limits $y_l^{av}$ and $T^{av}$ implies the frame-average expectations also have well defined limits,
with $\overline{y}_l = y_l^{av}$ and $\overline{T} = T^{av}$.  This holds, for example, whenever $y_l[r]$ and $T[r]$ are deterministically
bounded by finite constants, or when more general conditions hold that allow the Lebesgue dominated convergence theorem
to be applied \cite{williams-martingale}. 
Then the time average penalty per unit time associated with $y_l[r]$ (sampled only at renewal times for simplicity) 
satisfies with probability 1: 
\begin{eqnarray*}
\lim_{R\rightarrow\infty} \frac{\sum_{r=0}^{R-1}y_l[r]}{\sum_{r=0}^{R-1} T[r]}  =  \lim_{R\rightarrow\infty} \frac{\frac{1}{R}\sum_{r=0}^{R-1}y_l[r]}{\frac{1}{R}\sum_{r=0}^{R-1} T[r]} =  \frac{\overline{y}_l}{\overline{T}} 
\end{eqnarray*}
Therefore, the value $\overline{y}_l/\overline{T}$ indeed represents the limiting penalty per unit time associated with the process
$y_l[r]$. 

The problem (\ref{eq:p1})-(\ref{eq:p3}) seeks only to minimize a time average subject to time average constraints. 
The second problem we consider, more general than the first, seeks to maximize a \emph{concave and entrywise non-decreasing function} 
 $\phi(\bv{\gamma})$ of the time average
attribute vector ratio $\overline{\bv{x}}/\overline{T}$, where $\overline{\bv{x}} = (\overline{x}_1, \ldots, \overline{x}_M)$: 
\begin{eqnarray}
\mbox{Maximize:} & \phi(\overline{\bv{x}}/\overline{T}) \label{eq:phi1} \\
\mbox{Subject to:} & \overline{y}_l/\overline{T} \leq c_l \: \: \forall l \in \{1, \ldots, L\} \label{eq:phi2}  \\
& \pi[r] \in \script{P} \: \: \forall r \in \{0, 1, 2, \ldots\} \label{eq:phi3} 
\end{eqnarray}
where $\phi(\bv{\gamma})$ is a given concave and entrywise non-decreasing utility 
function defined over $\bv{\gamma} = (\gamma_1, \ldots, \gamma_M) 
\in \mathbb{R}^M$.  

We note that 
for some problems, one may be more interested in optimizing the per-frame average $\overline{y}_0$, rather than the time
average $\overline{y}_0/\overline{T}$, and such problems are treated in Section \ref{section:alt-form}. 

\subsection{Boundedness Assumptions} \label{section:boundedness} 

We assume $x_m[r]$, $T[r]$, and $y_0[r]$ have bounded conditional expectations, regardless of the policy. That is, there
are finite constants $x_{m}^{min}$, $x_m^{max}$, $T^{min}$, $T^{max}$, $y_0^{min}$, $y_0^{max}$ such that
for all $\pi[r] \in \script{P}$ and all $m \in \{1, \ldots, M\}$ we have: 
\begin{eqnarray*}
y_0^{min} \leq \expect{\hat{y}_0(\pi[r])|\pi[r]} \leq y_0^{max} \\
0 < T^{min} \leq \expect{\hat{T}(\pi[r])|\pi[r]} \leq T^{max} \\
x_m^{min} \leq \expect{\hat{x}_m(\pi[r])|\pi[r]} \leq x_m^{max}
\end{eqnarray*}
Define $\gamma_{m}^{min}$ and $\gamma_m^{max}$ by: 
\begin{eqnarray*}
 \gamma_m^{min} &\defequiv& \min[x_m^{min}/T^{min}, x_m^{min}/T^{max}] \\
 \gamma_m^{max} &\defequiv& \max[x_m^{max}/T^{max}, x_m^{max}/T^{max}]
 \end{eqnarray*}
 Define the hyper-rectangle $\script{R}$ by: 
 \begin{equation} \label{eq:R}
  \script{R} \defequiv \{ \bv{\gamma} \in \mathbb{R}^M | \gamma_m^{min} \leq \gamma_m \leq \gamma_m^{max} \: \: \forall m \in \{1, \ldots, M\} \} 
  \end{equation} 
 Then for any algorithm that chooses policies $\pi[r] \in \script{P}$ for all frames $r$,  it is not difficult to show that
 $\overline{x}_m[R]/\overline{T}[R] \in \script{R}$ for all $R \in \{1, 2, 3, \ldots\}$, 
 where $\overline{T}[R]$, $\overline{x}_m[R]$, $\overline{T}[R]$ are frame average expectations over the first $R$ frames, as defined
 by (\ref{eq:yav-r-frame}).
 
Finally, we assume the conditional second moments of $T[r]$, $x_m[r]$, and $y_l[r]$  (for $l \neq 0$) 
are finite, regardless of the policy. 
That is, there is a finite constant $\sigma_1$ such that for all $\pi[r] \in \script{P}$: 
\begin{eqnarray*}
\expect{\hat{T}(\pi[r])^2 | \pi[r]} &\leq& \sigma_1 \\ 
\expect{\hat{y}_l(\pi[r])^2 | \pi[r]} &\leq& \sigma_1 \: \: \forall l \in \{1, \ldots, L\} \\
\expect{\hat{x}_m(\pi[r])^2|\pi[r]} &\leq& \sigma_1 \: \: \forall m \in \{1, \ldots, M\} 
\end{eqnarray*}
%For technical reasons, we do not require the second moments of $y_0[r]$ to be finite. 

\subsection{Optimality of i.i.d. Algorithms}

We now state the problem (\ref{eq:p1})-(\ref{eq:p3}) more precisely, using $\limsup$s which do not require existence of a 
well defined limit: 
\begin{eqnarray} 
\mbox{Minimize:} & \limsup_{R\rightarrow\infty} \frac{\overline{y}_0[R]}{\overline{T}[R]} \label{eq:ls1} \\
\mbox{Subject to:} & \limsup_{R\rightarrow\infty} \frac{\overline{y}_l[R]}{\overline{T}[R]} \leq c_l \: \: \forall l \in \{1, \ldots, L\} \label{eq:ls2} \\
& \pi[r] \in \script{P} \: \: \forall r \in \{0, 1, 2, \ldots\} \label{eq:ls3} 
\end{eqnarray} 
Assume that the constraints (\ref{eq:ls2})-(\ref{eq:ls3}) are feasible, and 
define $ratio^{opt}$ as the infimum ratio in (\ref{eq:ls1}) over all algorithms that satisfy these constraints. 

Define an \emph{i.i.d. algorithm} as one that, at the beginning of each new frame $r \in \{0, 1, 2, \ldots\}$, chooses a policy 
$\pi[r]$ by independently and probabilistically selecting $\pi \in \script{P}$ according to some distribution that is the same
for all frames $r$.  Let $\pi^*[r]$ represent such an i.i.d. algorithm. Then the random variables 
$\{\hat{T}(\pi^*[r])\}_{r=0}^{\infty}$ are
independent and identically distributed (i.i.d.) over frames, as are $\{\hat{y}_l(\pi^*[r])\}_{r=0}^{\infty}$.  Thus,
by the law of large numbers, these have well defined time averages $\overline{T}^*$ and $\overline{y}_l^*$ with probability 1, 
where the averages are equal to the expectations over one frame. 

\begin{lem} \label{lem:optimality-iid} (Optimality over i.i.d. algorithms) If the constraints (\ref{eq:ls2})-(\ref{eq:ls3}) are feasible, 
then for any $\delta>0$, there exists an i.i.d. algorithm $\pi^*[r]$ that satisfies: 
\begin{eqnarray} 
\expect{\hat{y}_0(\pi^*[r])} \leq \expect{\hat{T}(\pi^*[r])}(ratio^{opt} + \delta) \label{eq:iid1} \\
\expect{\hat{y}_l(\pi^*[r])} \leq \expect{\hat{T}(\pi^*[r])}(c_l + \delta)  \: \: \forall l \in \{1, \ldots, L\} \label{eq:iid2} 
\end{eqnarray} 
\end{lem} 
\begin{proof}
See Appendix A.
\end{proof} 

\section{Optimizing Time Averages}

Here we develop an algorithm to treat the problem (\ref{eq:p1})-(\ref{eq:p3}).  To treat the constraints $\overline{y}_l/\overline{T} \leq c_l$, 
which are equivalent to the constraints $\overline{y}_l \leq c_l\overline{T}$, 
we define \emph{virtual queues} $Z_l[r]$ for $l \in \{1, \ldots, L\}$, with finite initial condition and with update equation: 
\begin{equation} \label{eq:z-update} 
Z_l[r+1] = \max[Z_l[r] + y_l[r] - c_lT[r], 0] \: \: \forall l \in \{1, \ldots, L\} 
\end{equation} 
The intuition is that if we can \emph{stabilize} the queue $Z_l[r]$, then the time average of the ``service process'' $c_l T[r]$ is
greater than or equal to the time average of the ``arrival process'' $y_l[r]$ (see also \cite{neely-energy-it} for application to 
\emph{virtual power queues} for meeting time average power constraints). 

Let $\bv{Z}[r] = (Z_1[r], \ldots, Z_L[r])$ be the vector of virtual queues, and 
define the following \emph{quadratic Lyapunov function} $L(\bv{Z}[r])$:
\begin{eqnarray*}
&L(\bv{Z}[r]) \defequiv \frac{1}{2}\sum_{l=1}^LZ_l[r]^2& 
\end{eqnarray*}
The value $L(\bv{Z}[r])$ is a scalar measure of the size of the queue backlogs.  The intuition is that if we 
can take actions that consistently push this value down, then queues can be stabilized.  
 Define the 
 \emph{frame-based conditional Lyapunov drift} $\Delta(\bv{Z}[r])$ by:
 \[ \Delta(\bv{Z}[r]) \defequiv \expect{L(\bv{Z}[r+1]) - L(\bv{Z}[r]) | \bv{Z}[r]} \]

\begin{lem} \label{lem:drift-compute} Under any control decision for choosing $\pi[r] \in \script{P}$, we have for all $r$
and all possible $\bv{Z}[r]$: 
\begin{eqnarray}
&\Delta(\bv{Z}[r]) \leq B + \expect{\sum_{l=1}^LZ_l[r][y_l[r] - c_lT[r]]|\bv{Z}[r]}& \label{eq:drift} 
\end{eqnarray}
where $B$ is a constant that satisfies for all $r$ and all possible $\bv{Z}[r]$: 
\begin{equation} \label{eq:B} 
B \geq \frac{1}{2}\sum_{l=1}^L\expect{(y_l[r] - c_lT[r])^2|\bv{Z}[r]} 
\end{equation} 
Such a constant $B$ exists by the boundedness assumptions in Section \ref{section:boundedness}. 
\end{lem} 
\begin{proof} 
Squaring (\ref{eq:z-update}) yields: 
\begin{eqnarray*}
Z_l[r+1]^2 &\leq& (Z_l[r] + y_l[r] - c_lT[r])^2 \\
&=& Z_l[r]^2 + (y_l[r] - c_lT[r])^2  \\
&& + 2Z_l[r](y_l[r] - c_lT[r])
\end{eqnarray*}
Taking conditional expectations, dividing by $2$,  
and summing over $l \in \{1, \ldots, L\}$ yields the result. 
\end{proof} 

\subsection{The Drift-Plus-Penalty Ratio Algorithm} \label{section:time-av-alg} 

Our \emph{Drift-Plus-Penalty Ratio Algorithm} is designed to minimize a sum of the variables on the right-hand-side 
of the drift bound (\ref{eq:drift}) and a penalty term, divided by an expected frame size, as in \cite{chihping-utility-round-robin}.  
The penalty term uses a non-negative
constant $V$ that will be shown to affect a performance tradeoff:
\begin{itemize}
\item (Policy Selection)  Every frame $r \in \{0, 1,2, \ldots\}$, observe the virtual queues $\bv{Z}[r]$ and choose a policy $\pi[r] \in \script{P}$ to minimize the following expression: 
\begin{eqnarray}
 \frac{\expect{ V\hat{y}_0(\pi[r]) + \sum_{l=1}^LZ_l[r]\hat{y}_l(\pi[r])|\bv{Z}[r]}}{\expect{\hat{T}(\pi[r])|\bv{Z}[r]}} \label{eq:dppr} 
 \end{eqnarray}
\item (Queue Update) Observe the resulting $\bv{y}[r]$ and $T[r]$ values, and 
update virtual queues $Z_l[r]$ by (\ref{eq:z-update}). 
\end{itemize}

Details on minimizing (\ref{eq:dppr}) are given in Section \ref{section:howto}. 
Rather than assuming we achieve the exact infimum of (\ref{eq:dppr}) over all policies $\pi[r] \in \script{P}$,  
it is useful to 
allow our decisions to come within 
an additive constant $C$ of the infimum. 

\begin{defn} A policy $\pi[r]$ is a \emph{$C$-additive approximation} for the problem (\ref{eq:dppr})  if for a given constant
$C\geq 0$ we
have: 
\begin{eqnarray*}
&&\frac{\expect{V\hat{y}_0(\pi[r]) + \sum_{l=1}^LZ_l[r]\hat{y}_l(\pi[r])|\bv{Z}[r]}}{\expect{\hat{T}(\pi[r])|\bv{Z}[r]}} \leq \\
&&C  + \inf_{\pi\in\script{P}}\left[\frac{\expect{ V\hat{y}_0(\pi) + \sum_{l=1}^LZ_l[r]\hat{y}_l(\pi)|\bv{Z}[r]}}{\expect{\hat{T}(\pi)|\bv{Z}[r]}} \right] 
\end{eqnarray*}
\end{defn} 

In Section \ref{section:pure-policies} it is shown that the infimum of (\ref{eq:dppr}) over $\pi\in \script{P}$ is the same as the infimum
over the extended class of probabilistically mixed strategies that choose a \emph{random} $\pi \in \script{P}$ according to some
distribution (exactly what  i.i.d. policies do every frame). Thus, if 
policy $\pi[r]$ is a $C$-additive approximation, then:
\begin{eqnarray}
&\expect{V\hat{y}_0(\pi[r]) + \sum_{l=1}^LZ_l[r]\hat{y}_l(\pi[r])|\bv{Z}[r]} \leq \nonumber \\
&\hspace{-.2in}\expect{\hat{T}(\pi[r])|\bv{Z}[r]} \left[ C +  \frac{\expect{V\hat{y}_0(\pi^*[r]) + \sum_{l=1}^LZ_l[r]\hat{y}_l(\pi^*[r])}}{\expect{\hat{T}(\pi^*[r])}}\right] \label{eq:ratio-good} 
\end{eqnarray} 
where $\pi^*[r]$ is any i.i.d. algorithm. Note that conditional expectations given $\bv{Z}[r]$ are the same as unconditional 
expectations under i.i.d. algorithms, because their decisions are independent of system history.

\begin{thm} \label{thm:alg1} (Algorithm Performance) Assume the constraints of 
problem (\ref{eq:ls1})-(\ref{eq:ls3}) are feasible. Fix constants $C\geq 0$, $V\geq 0$, and assume the 
above algorithm is implemented using any $C$-additive approximation every frame $r$
for the minimization in (\ref{eq:dppr}). Assume
initial conditions satisfy $\expect{L(\bv{Z}[0])} < \infty$. Then: 

a) For all $l \in \{1, \ldots, L\}$ we have: 
\begin{eqnarray}
 \limsup_{R\rightarrow\infty} \overline{y}_l[R]/\overline{T}[R] &\leq& c_l \: \: \forall l \in \{1, \ldots, L\} \label{eq:thm-mrs} \\
 \limsup_{R\rightarrow\infty} \frac{\sum_{r=0}^{R-1} y_l[r]}{\sum_{r=0}^{R-1}T[r]} &\leq& c_l \: \: \: (w.p.1) \label{eq:thm-wp1} 
\end{eqnarray}
where ``w.p.1'' stands for ``with probability 1.''

b) For all integers $R>0$ we have: 
\begin{eqnarray*}
\frac{\overline{y}_0[R]}{\overline{T}[R]} \leq ratio^{opt} + \frac{(B/\overline{T}[R]+C)}{V} + \frac{\expect{L(\bv{Z}[0])}}{VR\overline{T}[R]}
\end{eqnarray*}
and hence:
\begin{equation} \label{eq:thm-y0} 
 \limsup_{R\rightarrow\infty} \overline{y}_0[R]/\overline{T}[R] \leq ratio^{opt} + (B/T^{min} +C)/V 
 \end{equation} 
where $B$ is defined in (\ref{eq:B}), and $ratio^{opt}$ is the optimal solution to (\ref{eq:ls1})-(\ref{eq:ls3}). 
\end{thm} 

Thus, the algorithm satisfies all constraints, and the 
value of $V$ can be chosen appropriately large to make $(B/T^{min}+C)/V$ arbitrarily small, ensuring that 
the time average penalty is arbitrarily close to its optimal value $ratio^{opt}$.  The tradeoff in choosing a large value of 
$V$ comes in the size of the $Z_l[r]$ queues and the number of frames 
required for $\expect{Z_l[R]}/R$ to approach zero (which affects convergence time of the algorithm, see (\ref{eq:thm1-see}) in 
the proof).  In particular, in Appendix B it is shown that there are constants $F_1, F_2$ such that for all $l \in \{1, \ldots, L\}$ we have: 
\begin{equation} \label{eq:obtained-bound} 
  \frac{\overline{y}_l[R]}{\overline{T}[R]} \leq c_l + \frac{1}{T_{min}} \sqrt{\frac{F_1 + VF_2}{R} + \frac{\sum_{l=1}^L\expect{Z_l[0]^2}}{R^2}} 
  \end{equation} 
It is clear that the second term on the right-hand-side above vanishes as $R\rightarrow\infty$, but the number of frames
required for it to be small depends on the $V$ parameter.   
This bound holds for general problems.  A tighter 
bound can be obtained for problems with special structure (see (\ref{eq:see-structure}) in Appendix G). 

\begin{proof} (Theorem \ref{thm:alg1}) Consider any frame $r \in \{0, 1, 2, \ldots\}$. 
From (\ref{eq:drift}) we have: 
\begin{eqnarray*}
&\Delta(\bv{Z}[r]) + V\expect{\hat{y}_0(\pi[r])|\bv{Z}[r]} \leq B \\
&+ \expect{V\hat{y}_0(\pi[r]) + \sum_{l=1}^LZ_l[r][\hat{y}_l(\pi[r]) - c_l\hat{T}(\pi[r])]|\bv{Z}[r]} 
\end{eqnarray*}
Substituting (\ref{eq:ratio-good}) into the above yields: 
\begin{eqnarray}
&\Delta(\bv{Z}[r]) + V\expect{\hat{y}_0(\pi[r])|\bv{Z}[r]} \leq B +   \nonumber  \\
&\hspace{-.2in}\expect{\hat{T}(\pi[r])|\bv{Z}[r]}
 \left[ C + \frac{\expect{V\hat{y}_0(\pi^*[r]) + \sum_{l=1}^LZ_l[r]\hat{y}_l(\pi^*[r])}}{\expect{\hat{T}(\pi^*[r])}}\right] \nonumber \\
 &- \sum_{l=1}^LZ_l[r]c_l\expect{\hat{T}(\pi[r])|\bv{Z}[r]} \label{eq:plug-alg1} 
\end{eqnarray}
In the above inequality, $\pi[r]$ represents the $C$-additive approximate decision actually made, and $\pi^*[r]$ is from
any alternative i.i.d. algorithm.  Fixing any $\delta>0$, plugging 
the i.i.d. algorithm $\pi^*[r]$ from (\ref{eq:iid1})-(\ref{eq:iid2}) into the right-hand-side
of (\ref{eq:plug-alg1}), and letting $\delta \rightarrow 0$ yields: 
\begin{eqnarray}
&\Delta(\bv{Z}[r]) + V\expect{\hat{y}_0(\pi[r])|\bv{Z}[r]} \leq B  \nonumber \\
&+\expect{\hat{T}(\pi[r])|\bv{Z}[r]}[C + Vratio^{opt} ] \label{eq:ratio-good2proof} 
\end{eqnarray}
Taking expectations of the above yields: 
\begin{eqnarray}
\expect{L(\bv{Z}[r+1])} - \expect{L(\bv{Z}[r])}  + V\expect{\hat{y}_0(\pi[r])} \leq \nonumber \\
B + \expect{\hat{T}(\pi[r])}[ C+ Vratio^{opt}]  \label{eq:ratio-good3} 
\end{eqnarray}
Summing the above over $r \in \{0, \ldots, R-1\}$ for some integer $R>0$ and dividing by $R$ yields: 
\begin{eqnarray}
\frac{\expect{L(\bv{Z}[R])} - \expect{L(\bv{Z}[0])}}{R}  + V\overline{y}_0[R] \leq \nonumber \\
B +\overline{T}[R][C + Vratio^{opt}]  \label{eq:ratio-good4} 
\end{eqnarray}
Rearranging terms in the above and using the fact that $\expect{L(\bv{Z}[R])} \geq 0$ yields the result of part (b). 

To prove part (a), from (\ref{eq:ratio-good2proof}) 
there is a constant $F$ such that: 
\begin{eqnarray} 
\Delta(\bv{Z}[r]) \leq F \label{eq:thm1-see2} 
\end{eqnarray}
%where the constant $F$ is defined: 
%\begin{eqnarray*}
%F \defequiv B + CT^{max} + V[T^{max} ratio^{opt} - y_0^{min}] 
%\end{eqnarray*}
Thus, the drift of a quadratic Lyapunov function is bounded by a constant. Further, the second moments of 
per-frame changes in $Z_l[r]$ are bounded because of  
the second moment  assumptions on $y_l[r]$ and $T[r]$. 
It follows that (see \cite{lyap-opt2}): 
\begin{eqnarray}
\lim_{R\rightarrow\infty} \expect{Z_l[R]}/R &=& 0 \label{eq:mrs} \\
\lim_{R\rightarrow\infty} Z_l[R]/R &=& 0 \: \: \: (w.p.1) \label{eq:wp1} 
\end{eqnarray}
Now from the queue
update (\ref{eq:z-update}) we have for any frame $r$: 
\[ Z_l[r+1] \geq Z_l[r] + y_l[r] - c_lT[r] \]
Summing the above over $r \in \{0, \ldots, R-1\}$ for some integer $R>0$ 
yields:  
\begin{eqnarray} 
&Z_l[R] - Z_l[0] \geq \sum_{r=0}^{R-1} [y_l[r] - c_lT[r]]& \label{eq:det-bound} 
\end{eqnarray}
Taking expectations, dividing by $R$, and using $\expect{Z_l[0]} \geq 0$ yields for all integers $R>0$: 
\begin{equation*} 
\frac{\expect{Z_l[R]}}{R} \geq \overline{y}_l[R] - c_l\overline{T}[R] 
\end{equation*} 
Thus: 
\begin{equation} \label{eq:thm1-see}   
\frac{\overline{y}_l[R]}{\overline{T}[R]} \leq c_l + \frac{\expect{Z_l[R]}}{R\overline{T}[R]} \leq c_l + \frac{\expect{Z_l[R]}}{RT^{min}}
\end{equation} 
Taking limits of the above and using (\ref{eq:mrs}) proves (\ref{eq:thm-mrs}). A similar argument uses (\ref{eq:wp1}) to prove
(\ref{eq:thm-wp1}). 
\end{proof} 

Under a mild ``Slater-type'' assumption that ensures the constraints (\ref{eq:ls2}) are achievable with
 ``$\epsilon$-slackness,''  the queues $Z_l[R]$ can be shown to be \emph{strongly stable}, in the sense that the
 time average expectation is 
 bounded by $O(V)$.  
 If further mild \emph{fourth moment boundedness assumptions} hold for $y_l[r]$ and $T[r]$ then 
 the same bound (\ref{eq:thm-y0}) can be shown to hold for pure time averages
  with probability 1 (see \cite{lyap-opt2} and Appendix H). 
   
 \section{Utility Optimization} 
 
 Consider now the problem (\ref{eq:phi1})-(\ref{eq:phi3}), which seeks to maximize $\phi(\overline{\bv{x}}/\overline{T})$ 
 subject to $\overline{y}_l/\overline{T} \leq c_l$ for all $l \in \{1, \ldots, L\}$.   We 
 transform this problem of maximizing a function of a time average ratio into a problem of
 the type (\ref{eq:p1})-(\ref{eq:p3}).  The following variation on Jensen's inequality is crucial in this transformation: 
 
 \begin{lem} \label{lem:jensen} (Variation on Jensen's Inequality) 
 Let $\phi(\bv{\gamma})$ be any continuous and concave function defined over $\bv{\gamma} \in \script{R}$
 for some closed and bounded hyper-rectangle $\script{R}$. 
 
 a) Let $(T, \bv{\gamma})$ be any random vector (with an arbitrary joint distribution) that satisfies 
 $T>0$ and $\bv{\gamma} \in \script{R}$ with probability 1.  Assume that  $0 < \expect{T} < \infty$. 
 Then: 
 \[  \frac{\expect{T\phi(\bv{\gamma})}}{\expect{T}} \leq \phi\left(\frac{\expect{T\bv{\gamma}}}{\expect{T}}\right) \]
 
 b) Let $(T[r], \bv{\gamma}[r])$ be a sequence of arbitrarily correlated 
  random vectors for $r \in \{0, 1, 2, \ldots\}$. Assume that $T[r]>0$, $\bv{\gamma}[r] \in \script{R}$ for all $r$ (with probability 1), 
  and: 
  \[ 0 < T^{min} \leq \expect{T[r]} \leq T^{max}  < \infty \: \: \forall r \in \{0, 1, 2, \ldots\} \]
    Then for any $R>0$: 
  \[ \frac{\frac{1}{R}\sum_{r=0}^{R-1}\expect{T[r]\phi(\bv{\gamma}[r])}}{\frac{1}{R}\sum_{r=0}^{R-1}\expect{T[r]}} \leq 
  \phi\left(\frac{\frac{1}{R}\sum_{r=0}^{R-1}\expect{T[r]\bv{\gamma}[r]}}{\frac{1}{R}\sum_{r=0}^{R-1}\expect{T[r]}}\right) \]
  Furthermore, assuming that the limits $\overline{T\phi(\bv{\gamma})}$ and $\overline{T\bv{\gamma}}$
   defined below exist, we have: 
  \begin{equation} \label{eq:crucial} 
   \overline{T\phi(\bv{\gamma})}/\overline{T} \leq \phi(\overline{T\bv{\gamma}}/\overline{T}) 
   \end{equation} 
    where: 
 \begin{eqnarray*}
\overline{T\phi(\bv{\gamma})} &\defequiv& \lim_{R\rightarrow\infty}  \frac{1}{R}\sum_{r=0}^{R-1}\expect{T[r]\phi(\bv{\gamma}[r])} \\
\overline{T\bv{\gamma}} &\defequiv& \lim_{R\rightarrow\infty}  \frac{1}{R}\sum_{r=0}^{R-1}\expect{T[r]\bv{\gamma}[r]} 
\end{eqnarray*}
 \end{lem} 
 \begin{proof} 
 Part (b) follows immediately from part (a) by defining the random vector $(T, \bv{\gamma})$ to be 
 $(T[J], \bv{\gamma}[J])$, where $J$ is a uniformly distributed integer in $\{0, 1, \ldots, R-1\}$ that is independent
 of the $(T[r], \bv{\gamma}[r])$ process.  Part (a) is proven in Appendix E. 
 \end{proof}

 Now define an auxiliary vector $\bv{\gamma}[r]= (\gamma_1[r], \ldots, \gamma_M[r])$, to be chosen in the set $\script{R}$ defined in (\ref{eq:R}) on every 
 frame $r$.  
 \begin{lem} (Equivalent Transformation) The problem (\ref{eq:phi1})-(\ref{eq:phi3}) is equivalent to the following transformed problem: 
 \begin{eqnarray}
 \mbox{Maximize:} & \overline{T\phi(\bv{\gamma})}/\overline{T} \label{eq:t1} \\
 \mbox{Subject to:} & \overline{x}_m \geq \overline{T\gamma_m} \: \: \forall m \in \{1, \ldots, M\} \label{eq:t2} \\
 & \overline{y}_l/\overline{T} \leq c_l \: \: \forall l \in \{1, \ldots, L\} \label{eq:t3} \\
 & \bv{\gamma}[r] \in \script{R} \: \: \forall r \in \{0, 1, 2, \ldots\} \label{eq:t4} \\
 & \pi[r] \in \script{P} \: \: \forall r \in \{0, 1, 2, \ldots\} \label{eq:t5} 
 \end{eqnarray}
 \end{lem} 
 \begin{proof} We briefly sketch the proof:  Let $\pi^*[r]$, $\bv{\gamma}^*[r]$ be a policy that optimally solves the above transformed problem, 
 and assume for simplicity it yields well defined time averages $\overline{T}^*$, $\overline{y}_l^*$, $\overline{x}_m^*$, 
 $\overline{T^*\phi(\bv{\gamma}^*)}$, $\overline{T^*\bv{\gamma}^*}$, and  
 optimal utility $util^* = \overline{T^*\phi(\bv{\gamma}^*)}/\overline{T}^*$.  Then the policy $\pi^*[r]$ also satisfies all constraints of problem
 (\ref{eq:phi1})-(\ref{eq:phi3}), and yields: 
 \[ \phi(\overline{\bv{x}}^*/\overline{T}^*) \geq \phi(\overline{T^*\bv{\gamma^*}}/\overline{T^*}) \geq \overline{T^*\phi(\bv{\gamma}^*)}/\overline{T}^* \defequiv
 util^* \]
 where the first inequality above
 holds by (\ref{eq:t2}) and the entrywise non-decreasing property of $\phi(\bv{\gamma})$, and the second holds
 by (\ref{eq:crucial}).  Thus, the optimal utility of problem (\ref{eq:phi1})-(\ref{eq:phi3}) is greater than or equal to that of the 
 transformed problem.   A similar argument shows it is also less than or equal to the optimal utility of the transformed problem.  
 \end{proof} 
 
 The transformed problem (\ref{eq:t1})-(\ref{eq:t5}) has the structure of the problem (\ref{eq:p1})-(\ref{eq:p3}) if we define
 $y_0[r] \defequiv -T[r]\phi(\bv{\gamma}[r])$, write the constraints (\ref{eq:t2}) as $\overline{T\gamma_m - x_m} \leq 0$, 
 and define policy decision $\pi'[r] \defequiv (\pi[r], \bv{\gamma}[r]) \in \script{P} \times \script{R}$. The resulting algorithm is thus the same as that given 
 in Section \ref{section:time-av-alg}, and for this context it is given as follows:  For the constraints (\ref{eq:t3}), use the same
 virtual queues $Z_l[r]$ defined in (\ref{eq:z-update}). For the constraints (\ref{eq:t2}), define virtual queues
$G_m[r]$ for $m \in \{1, \ldots, M\}$ by: 
\begin{equation} \label{eq:g-update} 
G_m[r+1] = \max[G_m[r] + T[r]\gamma_m[r] - x_m[r], 0] 
\end{equation} 
Define $\bv{G}[r] \defequiv (G_1[r], \ldots, G_M[r])$. 
The drift-plus-penalty ratio to minimize every frame $r$ is then: 
\begin{eqnarray*}
 \frac{\expect{-V\hat{T}(\pi[r])\phi(\bv{\gamma}[r]) + \sum_{l=1}^LZ_l[r]\hat{y}_l(\pi[r]) | \bv{Z}[r]}}{\expect{\hat{T}(\pi[r])|\bv{Z}[r]}} \\
 + \frac{\expect{\sum_{m=1}^MG_m[r][\hat{T}(\pi[r])\gamma_m[r] - \hat{x}_m(\pi[r])]|\bv{Z}[r]}}{\expect{\hat{T}(\pi[r])|\bv{Z}[r]}} 
 \end{eqnarray*}
 It is easy to see that the above can be minimized by separately choosing $\bv{\gamma}[r] \in \script{R}$ and $\pi[r] \in \script{P}$
 to minimize their respective terms, and that $\hat{T}(\pi[r])$ cancels out of the auxiliary variable decisions.  The resulting algorithm is thus
 to observe $\bv{Z}[r]$ and $\bv{G}[r]$ every frame $r\in \{0, 1,2, \ldots\}$ and perform the following: 
 \begin{itemize} 
 \item (Auxiliary Variables) Choose $\bv{\gamma}[r] \in \script{R}$ to maximize: 
 \begin{eqnarray*}
 &V\phi(\bv{\gamma}[r]) - \sum_{m=1}^MG_m[r]\gamma_m[r]&
 \end{eqnarray*}
 \item (Policy Selection) Choose $\pi[r] \in \script{P}$ to minimize: 
 \[ \frac{\expect{ \sum_{l=1}^LZ_l[r]\hat{y}_l(\pi[r]) - \sum_{m=1}^MG_m[r]\hat{x}_m(\pi[r])|\bv{Z}[r]}}{\expect{\hat{T}(\pi[r])|\bv{Z}[r]}} \]
 \item (Virtual Queue Update) Update $\bv{Z}[r]$ by (\ref{eq:z-update}) and $\bv{G}[r]$ by (\ref{eq:g-update}). 
 \end{itemize} 
 The auxiliary variable update is a simple deterministic
 maximization of a concave function over a hyper-rectangle, and can be separated
 into $M$ optimizations of single-variable concave functions over an interval 
 if the utility function has the form
 $\phi(\bv{\gamma}) = \sum_{m=1}^M\phi_m(\gamma_m)$. 
 The policy selection step is again an optimization of a ratio of expectations and can be done as described in Section 
 \ref{section:howto}. 
 
 We define a $C$-additive approximation of the above algorithm as one that, every frame $r$, chooses $\pi[r] \in \script{P}$ 
 to yield an expectation of ratios in the policy selection step that is within $C$ of the infimum.
 To explicitly describe the performance of the above algorithm, we write the problem (\ref{eq:phi1})-(\ref{eq:phi3}) more precisely
 using $\limsup$s: 
 \begin{eqnarray} 
 \mbox{Minimize:} & \limsup_{R\rightarrow\infty} \phi(\overline{\bv{x}}[R]/\overline{T}[R]) \label{eq:tran1} \\
 \mbox{Subject to:} & \limsup_{R\rightarrow\infty} \frac{\overline{y}_l[R]}{\overline{T}[R]} \leq c_l \: \: \forall l \in \{1, \ldots, L\} \label{eq:tran2} \\
 & \pi[r] \in \script{P} \: \: \forall r \in \{0, 1, 2, \ldots\} \label{eq:tran3} 
 \end{eqnarray} 
 Assuming the constraints of the above problem are feasible, define $util^{opt}$ as the supremum value of (\ref{eq:tran1}) 
 over all algorithms that satisfy (\ref{eq:tran2})-(\ref{eq:tran3}). 
 
 \begin{thm} \label{thm:utility} Suppose the constraints of problem (\ref{eq:tran1})-(\ref{eq:tran3}) are feasible, and 
 a $C$-additive approximation is used every frame  (for $C\geq0$).  Then: 
 
 a) For all $l \in \{1, \ldots, L\}$ we have: 
\begin{eqnarray*}
 \limsup_{R\rightarrow\infty} \overline{y}_l[R]/\overline{T}[R] &\leq& c_l \: \: \forall l \in \{1, \ldots, L\} \\
 \limsup_{R\rightarrow\infty} \frac{\sum_{r=0}^{R-1} y_l[r]}{\sum_{r=0}^{R-1}T[r]} &\leq& c_l \: \: \: (w.p.1)  
\end{eqnarray*}

b) The achieved utility satisfies: 
\[  \liminf_{R\rightarrow\infty} \phi\left(\frac{\overline{\bv{x}}[R]}{\overline{T}[R]}\right) \geq util^{opt} - \frac{D}{VT^{min}}  - \frac{C}{V} \]  
where $D$ is a constant that satisfies for all $r$ and all possible $\bv{Z}[r]$: 
 \begin{eqnarray}
  D &\geq&   \frac{1}{2}\sum_{l=1}^L\expect{[y_l[r] - c_lT[r]]^2|\bv{Z}[r], \bv{G}[r]} \nonumber  \\
 &&  + \frac{1}{2} \sum_{m=1}^M\expect{[T[r]\gamma_m[r] - x_m[r]]^2|\bv{Z}[r], \bv{G}[r]} \label{eq:D}
  \end{eqnarray}
  Such a constant $D$ exists by the boundedness assumptions in Section \ref{section:boundedness}. 
\end{thm} 
  \begin{proof} 
 See Appendix F. 
 \end{proof} 
 
\section{Optimizing the Ratio of Expectations} \label{section:howto} 

Here we show how to minimize the ratio of 
expectations given in (\ref{eq:dppr}) (and also in the policy selection stage of the previous
section). These problems can be written more generally 
as choosing a policy $\pi[r] \in \script{P}$ to minimize the ratio: 
\[ \frac{\expect{a(\pi)}}{\expect{b(\pi)}} \]
where $a(\pi), b(\pi)$ are random functions of $\pi \in \script{P}$.  
The function $b(\pi)$ is equal to $\hat{T}(\pi)$, and is strictly positive and satisfies the following
for all $\pi \in \script{P}$: 
\[ 0 < T^{min} \leq \expect{b(\pi)|\pi} \leq T^{max} < \infty \]
The function
$a(\pi)$ depends on $\bv{Z}[r]$, and 
the above expectations are implicitly conditioned on $\bv{Z}[r]$, although we suppress this notation for
simplicity.  
Define $\theta^*$ as the optimal ratio: 
\[ \theta^* \defequiv \inf_{\pi \in \script{P}} \left[\frac{\expect{a(\pi)}}{\expect{b(\pi)}}\right] \]
If the expectation $\expect{b(\pi)}$ is the same for all $\pi \in \script{P}$ (such as when the frame size is independent of the policy), 
then $\theta^*$ is obtained by infimizing 
the numerator $\expect{a(\pi)}$.  This is typically easier 
(often involving learning for 
\emph{stochastic shortest path} computations \cite{bertsekas-neural}\cite{neely-mdp-cdc09}).   Otherwise, the following simple lemma is useful.

\begin{lem} \label{lem:bisection} We have: 
\begin{equation} 
 \inf_{\pi\in\script{P}} \expect{a(\pi) - \theta^*b(\pi)} = 0 \label{eq:thetagood} 
\end{equation} 
Further, for any real number $\theta$, we have: 
\begin{eqnarray}
\inf_{\pi\in\script{P}} \expect{a(\pi) - \theta b(\pi)} < 0 & \mbox{ if $\theta > \theta^*$} \label{eq:thetaup} \\
\inf_{\pi\in\script{P}} \expect{a(\pi) - \theta b(\pi)} > 0 & \mbox{ if $\theta < \theta^*$} \label{eq:thetadown} 
\end{eqnarray}
\end{lem} 
\begin{proof} 
We first assume the result of (\ref{eq:thetagood}) and use it to 
prove (\ref{eq:thetaup})-(\ref{eq:thetadown}). 
Suppose
that $\theta > \theta^*$. We then have for any $\pi \in \script{P}$: 
\begin{eqnarray*}
 \expect{a(\pi) - \theta b(\pi)} &=&  \expect{a(\pi) - \theta^*b(\pi) - (\theta-\theta^*)b(\pi)} \\
 &\leq& \expect{a(\pi) - \theta^*b(\pi)} - (\theta - \theta^*)T^{min}
\end{eqnarray*}
Thus: 
\begin{eqnarray*}
 \inf_{\pi\in\script{P}} \expect{a(\pi) - \theta b(\pi)} &\leq& \inf_{\pi\in\script{P}}\expect{a(\pi) - \theta^*b(\pi)} \\
 && - (\theta - \theta^*)T^{min} \\
 &=& 0 - (\theta - \theta^*)T^{min} < 0
 \end{eqnarray*}
 where the equality holds by (\ref{eq:thetagood}). This proves (\ref{eq:thetaup}). 
 
 Now suppose that $\theta < \theta^*$. Then for any $\pi \in \script{P}$: 
 \begin{eqnarray*}
 \expect{a(\pi) - \theta b(\pi)} &=&  \expect{a(\pi) - \theta^*b(\pi) +  (\theta^*-\theta)b(\pi)} \\
 &\geq& \expect{a(\pi) - \theta^*b(\pi)} +  (\theta^* - \theta)T^{min}
\end{eqnarray*}
Taking infimums of both sides, again using (\ref{eq:thetagood}), proves: 
\begin{eqnarray*}
 \inf_{\pi\in\script{P}} \expect{a(\pi) - \theta b(\pi)} \geq 0 + (\theta^*-\theta)T^{min} > 0
 \end{eqnarray*}
 This proves (\ref{eq:thetadown}). 
 
 It remains only to prove (\ref{eq:thetagood}).  We have for any policy $\pi \in \script{P}$: 
 \[ \frac{\expect{a(\pi)}}{\expect{b(\pi)}} \geq \inf_{\pi\in\script{P}} \left[\frac{\expect{a(\pi)}}{\expect{b(\pi)}}\right] \defequiv \theta^* \]
 Therefore, because $\expect{b(\pi)}>0$, we have for all $\pi \in \script{P}$:  
 \[ \expect{a(\pi)} - \theta^*\expect{b(\pi)} \geq 0 \]
 and hence: 
 \[ \inf_{\pi\in\script{P}}\expect{a(\pi) - \theta^*b(\pi)} \geq 0 \]
 It remains only to prove the reverse inequality. Fix $\delta>0$.  By definition of $\theta^*$ as the infimum ratio, 
 there is a policy $\pi^*\in\script{P}$ that satisfies: 
 \[ \frac{\expect{a(\pi^*)}}{\expect{b(\pi^*)}} \leq \theta^* + \delta \]
 Thus:
 \[ \expect{a(\pi^*)} \leq  \theta^*\expect{b(\pi^*)} + \delta\expect{b(\pi^*)} \]
 and so: 
 \[ \expect{a(\pi^*)} - \theta^*\expect{b(\pi^*)} \leq \delta T^{max} \]
 Because $\pi^*$ is just a particular algorithm in $\script{P}$, it follows that:  
 \[ \inf_{\pi\in\script{P}} \expect{a(\pi) - \theta^*b(\pi)} \leq \delta T^{max} \]
 This holds for all $\delta>0$. Taking a limit as $\delta\rightarrow0$ completes the proof. 
\end{proof} 

\subsection{The Bisection Algorithm} 

Lemma \ref{lem:bisection} immediately leads to the following simple bisection algorithm:  Suppose we have upper
and lower bounds $\theta_{min}$ and $\theta_{max}$, so that we know $\theta_{min} \leq \theta^* \leq \theta_{max}$. 
Then we can define $\theta = (\theta_{min} + \theta_{max})/2$, and compute the value of $\inf_{\pi\in \script{P}} \expect{a(\pi) - \theta b(\pi)}$. 
If the result is $0$, then $\theta = \theta^*$.  If positive, then $\theta < \theta^*$, and otherwise $\theta > \theta^*$. We can then
refine our upper and lower bounds.  This leads to a simple iterative algorithm where the distance between the upper and lower
bounds decreases by a factor of 2 on each iteration. It thus approaches the optimal $\theta^*$ value exponentially fast.  Each step of the iteration involves minimizing an expectation, rather than a ratio of expectations. 

\subsection{Optimizing over Pure Policies} \label{section:pure-policies}

Note that for any set of policies $\script{S}$, 
Lemma \ref{lem:bisection} implies that $\inf_{\pi\in\script{S}} \expect{a(\pi) - \theta b(\pi)} = 0$ if and only if $\theta = \inf_{\pi\in\script{S}} \expect{a(\pi)}/\expect{b(\pi)}$. 
Now suppose we have a set of policies $\script{P}^{pure}$ that we call \emph{pure policies}, and that the policy space $\script{P}$ 
consists of 
all pure policies as well as all ``mixtures'' (or convex combinations) of pure policies, being policies that choose a pure policy in $\script{P}^{pure}$ with some
particular probability distribution.  More generally, define $\Omega$ as the set of all vectors $(\expect{a(\pi)}, \expect{b(\pi)})$ achievable over $\pi \in \script{P}^{pure}$, and suppose the set of all $(\expect{a(\pi)}, \expect{b(\pi)})$ achievable over $\pi \in \script{P}$ is equal to the convex hull of $\Omega$.  Recall that $\theta^*$ is the infimum 
ratio of $\expect{a(\pi)}/\expect{b(\pi)}$ over $\pi \in \script{P}$. 
Then: 
\begin{eqnarray*}
0 = \inf_{\pi\in\script{P}} \expect{a(\pi) - \theta^* b(\pi)} 
&=& \inf_{(a,b) \in Conv(\Omega)} [a - \theta^*b] \\
&=& \inf_{(a,b) \in \Omega} [a - \theta^* b] \\
&=& \inf_{\pi\in\script{P}^{pure}} \expect{a(\pi) - \theta^*b(\pi)} 
\end{eqnarray*}
where the third inequality holds because the infimum of a linear function over the convex hull of a
set is equal to the infimum over
the set itself. 
It follows 
that $\theta^*$ is also the infimum ratio of $\expect{a(\pi)}/\expect{b(\pi)}$ over $\pi \in \script{P}^{pure}$.  

This means that to achieve
the infimum ratio over policies $\pi \in \script{P}$, it suffices to restrict our search to pure policies. 
%In particular, if our policy space $\script{P}$ includes all mixtures of a finite set of pure policies $\script{P}^{pure} = \{\pi_1, \pi_2, \ldots, \pi_K\}$ (for some finite integer $K$), then we can simply evaluate the ratio $\expect{a(\pi)}/\expect{b(\pi)}$ over each 
%of the $K$ pure policies and select the minimizing one.   We do not need to test the ratio for any of the 
%(infinite number) of mixed policies. 

\subsection{Optimizing with Initial Information} \label{section:initial-info} 

Suppose at the beginning of each frame, we observe a vector $\bv{\eta}[r]$ of \emph{initial information} that can affect
the penalties and frame size. Suppose that $\{\bv{\eta}[r]\}_{r=0}^{\infty}$ is i.i.d. over frames.  Each policy $\pi \in \script{P}$ first
observes $\bv{\eta}[r]$ and then chooses a \emph{sub-policy} $\pi' \in \script{P}_{\bv{\eta}[r]}$, where $\script{P}_{\bv{\eta}[r]}$ is a space
that possibly depends on $\bv{\eta}[r]$.  To minimize $\expect{a(\pi)}$, it suffices to observe $\bv{\eta}[r]$ and choose $\pi' \in \script{P}_{\bv{\eta}[r]}$ to minimize the conditional expectation $\expect{a(\pi')| \bv{\eta}[r]}$.  However, this is not necessarily true
for minimizing the ratio $\expect{a(\pi)}/\expect{b(\pi)}$. 
%That is, minimizing the ratio of expectations is \emph{not} always achieved by minimizing the ratio of conditional 
%expectations given the observed $\bv{\eta}[r]$ (see \cite{sno-text} for a simple counter-example). 

A correct approach is the following:  If $\theta^*$ is known, we can simply choose $\pi' \in \script{P}_{\bv{\eta}[r]}$ to minimize: 
\[ \expect{a(\pi') - \theta^* b(\pi')|\bv{\eta}[r]} \]
If $\theta^*$ is unknown, we can carry out the bisection routine.  Let $\theta$ be the midpoint in the current iteration.  We must
compute: 
\begin{equation} \label{eq:iterated-expectations} 
\inf_{\pi \in \script{P}} \expect{a(\pi) - \theta b(\pi)} = \expect{\inf_{\pi'\in\script{P}_{\bv{\eta}[r]}}\expect{a(\pi') - \theta b(\pi')| \bv{\eta}[r]}} 
\end{equation} 
The infimizing decision $\pi'$ can be made by observing $\bv{\eta}[r]$,  without requiring knowledge of its probability
distribution.  However, the 
value in (\ref{eq:iterated-expectations}) cannot be computed without knowledge of this distribution.  Instead, 
suppose we have $W$ i.i.d. samples $\{\bv{\eta}_w\}_{w=1}^W$.  We can then approximate the value in (\ref{eq:iterated-expectations}) by the function $val(\theta)$ defined below: 
\begin{eqnarray}
val(\theta) \defequiv \frac{1}{W}\sum_{w=1}^W\inf_{\pi'\in\script{P}_{\bv{\eta}_w}}\expect{a(\pi') - \theta b(\pi')| \bv{\eta}_w} \label{eq:val-theta} 
\end{eqnarray}

By the law of large numbers, $val(\theta)$ approaches the exact value of (\ref{eq:iterated-expectations}) with a large choice of $W$.
The bisection routine can be carried out using the $val(\theta)$ approximation, being sure to use the same samples at each step of the iteration (but different samples on each frame $r$). 
Note that $val(\theta)$ is non-increasing in $\theta$, so the bisection will converge provided that it is initialized so that
$val(\theta_{min}) \geq 0$ and $val(\theta_{max}) \leq 0$.  If we cannot independently generate $W$ samples, we 
use the $W$ past observed values of $\bv{\eta}[r]$ from previous frames.  
There is a subtle issue here, as these past values have influenced system performance
and are thus correlated with the current 
$a(\pi)$ and $b(\pi)$ functions.  However, a \emph{delayed queue argument} similar to 
that given in \cite{neely-mwl-arxiv} shows these past values can still be used. 

\section{Alternative Algorithms without Ratio Minimization} \label{section:alts} 

We first present an alternative formulation to (\ref{eq:p1})-(\ref{eq:p3}) that is easier and does not require minimizing
a ratio of expectations every slot.  We then present an alternative algorithm to the original problem (\ref{eq:p1})-(\ref{eq:p3}) that
does not require a ratio of minimizations, but which yields a less explicit convergence result.

\subsection{Alternative Formulation} \label{section:alt-form}

Note that constraints of the form $\overline{y}_l \leq 0$ are equivalent to $\overline{y}_l/\overline{T} \leq c_l$ in the special case $c_l=0$,
and thus can be handled using the framework of this paper.  Thus, if $f[r]$ is some penalty, and if we desire the constraint
$\overline{f} \leq 6.7$, then we can define $y[r] \defequiv f[r] - 6.7$ and note that the desired constraint is equivalent to 
$\overline{y} \leq 0$.  In other words, constraints of the form 
$\overline{y}_l/\overline{T} \leq c_l$ are more general than constraints 
of the form $\overline{y}_l \leq 0$ or $\overline{y}_l \leq c$ for some constant $c$, and contain these as special cases.

Now consider the following problem structure: 
\begin{eqnarray*}
\mbox{Minimize:} & \overline{y}_0 \\
\mbox{Subject to:} & \overline{y}_l/\overline{T} \leq c_l \: \: \forall l \in \{1, \ldots, L\} \\
& \pi[r] \in \script{P} \: \: \forall r \in \{0, 1,2, \ldots \} 
\end{eqnarray*}
Such a problem has a different structure than the problem (\ref{eq:p1})-(\ref{eq:p3}), and is \emph{easier to solve} as it does not require a ratio of expectations. It can be solved using the same virtual queues $Z_l[r]$ in (\ref{eq:z-update}), but every frame
$r$ observing $\bv{Z}[r]$ and selecting a policy 
$\pi[r] \in \script{P}$ to minimize the following expression: 
\begin{eqnarray*}
&\mathbb{E}\{V\hat{y}_0(\pi[r]) + \sum_{l=1}^LZ_l[r][\hat{y}_l(\pi[r]) - c_l\hat{T}(\pi[r])]|\bv{Z}[r]\}&
\end{eqnarray*}
Analysis of this algorithm is given in Appendix C.  

\subsection{Alternative Algorithm} \label{section:alternative} 
The following is an alternative algorithm for the original problem (\ref{eq:p1})-(\ref{eq:p3}) that does not
require a ratio minimization (and hence does not require a bisection step):  Use the same virtual queues 
$Z_l[r]$ in (\ref{eq:z-update}).  Define $\theta[0]=0$, and define $\theta[R]$ for $R \in \{1, 2, 3, \ldots\}$ by: 
\begin{eqnarray} 
&\theta[R] \defequiv \sum_{r=0}^{R-1} y_0[r]/\sum_{r=0}^{R-1}T[r]& \label{eq:theta-timeav}
\end{eqnarray} 
Every frame
$r$, observe $\bv{Z}[r]$ and $\theta[r]$ and select a policy 
$\pi[r] \in \script{P}$ to minimize the following expression: 
\begin{eqnarray}
&\mathbb{E}\{V[\hat{y}_0(\pi[r]) - \theta[r]\hat{T}(\pi[r])]|\bv{Z}[r], \theta[r]\} \label{eq:alt-metric}  \\
 & +  \mathbb{E}\{\sum_{l=1}^LZ_l[r][\hat{y}_l(\pi[r]) - c_l\hat{T}(\pi[r])]|\bv{Z}[r], \theta[r]\} \nonumber
 \end{eqnarray}
It is shown in Appendix D 
that all constraints are met, and that 
if $\theta[r]$ converges to a constant with probability 1, then with probability 1: 
\begin{eqnarray*}
&\lim_{R\rightarrow\infty} \sum_{r=0}^{R-1} y_0[r]/\sum_{r=0}^{R-1}T[r] \leq ratio^{opt} + O(1/V)& 
 \end{eqnarray*}
The disadvantage is that the convergence time is not as clear as that given in part (b) of 
Theorem \ref{thm:alg1}.  
Further, use of the time average (\ref{eq:theta-timeav}) makes it difficult to adapt to
changes in system parameters, so that it may be better to approximate (\ref{eq:theta-timeav}) with a moving 
average or an exponentially decaying average. 

\section{Simulations for a Task Processing Network} \label{section:simulation}

\begin{figure}[cht]
   \centering
   \includegraphics[height=.6in, width=3in]{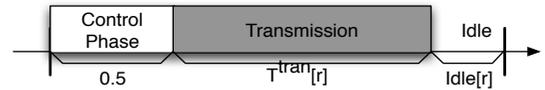} % requires the graphicx package
   \caption{An illustration of the 3 phases of a renewal frame $r \in \{0, 1,2, \ldots\}$.}
   \label{fig:task-structure}
\end{figure}

Here we provide a simple task processing example. An infinite sequence of
tasks must be processed one at a time with the help of a network of 5 wireless devices.  
This applies, for example, in scenarios similar to  \cite{network-corroboration-arxiv} where each 
new task represents an event that is sensed by the wireless devices (each at different sensing qualities \cite{bisdikian-qoi}), 
and we must select which device reports the event information. 
The renewal structure is shown in Fig. \ref{fig:task-structure}. At the beginning
of each new task $r$, a period of $0.5$ time units is expended to communicate control information about the
task.  Each of the 5 devices expends $0.5$ units of energy in this control phase. 
At the end of this phase, the network controller obtains a vector $\bv{\eta}[r]$ of parameters for task $r$. 
The vector $\bv{\eta}[r]$ has the form: 
\[\bv{\eta}[r] = [(qual_1[r], T_1^{tran}[r]), \ldots, (qual_5[r], T_5^{tran}[r])] \]
where for each $l \in \{1, \ldots, 5\}$, 
$qual_l[r]$ is a real number representing the \emph{information quality} if device $l$ is  chosen to process
task $r$, and $T_l^{tran}[r]$ is the \emph{transmission time} required for device $l$ to transmit the corresponding
information to a receiving station.   The controller must choose one of the 5 devices to process the task, and must also 
choose the amount of \emph{idle time} at the end of the frame (chosen within the interval $[0, I^{max}]$ for some constant
$I^{max}>0$), so that the 
policy decision $\pi[r]$ has the form: 
\[ \pi[r] = (l[r], Idle[r]) \in \{1, 2, 3, 4, 5\} \times \{I \in \mathbb{R} | 0 \leq I \leq I^{max}\} \] 
 
Define $P^{tran}$ as the power expenditure associated with wireless transmission.  
The chosen device $l[r]$ expends $P^{tran} \times T_{l[r]}^{tran}$ units of energy in the transmit phase, while all other
 devices $l \neq l[r]$ expend no energy in this phase.  None of the devices expend energy in the idle phase, which 
 helps to limit the average power expenditure in the system. 
 
 The goal is to maximize the \emph{quality of information (q.o.i) per unit time} subject to an average power constraint of $0.25$
 at each device. 
 Define $\hat{y}_0(\pi[r])$ as $-1$ times the q.o.i.  obtained 
 for task $r$, $\hat{y}_l(\pi[r])$ as the energy expended 
 by device $l$ on task $r$, and $\hat{T}(\pi[r])$ as the frame duration for task $r$: 
 \begin{eqnarray*}
 \hat{y}_0(\pi[r]) &\defequiv& -qual_{l[r]}[r] \\
 \hat{y}_l(\pi[r]) &\defequiv& 0.5 + P^{tran}T_l^{tran}[r]1_{\{l[r] = l\}} \: \: \forall l \in \{1, \ldots, 5\}\\
 \hat{T}(\pi[r]) &\defequiv& 0.5 + T_{l[r]}^{tran}[r] + Idle[r]
 \end{eqnarray*}
 where $1_{\{l[r] = l\}}$ is an indicator function that is $1$ if $l[r]=l$
 and $0$ else.    The problem is then to minimize $\overline{y}_0/\overline{T}$ subject to $\overline{y}_l/\overline{T} \leq 0.25$ for
 all  $l \in \{1, \ldots, 5\}$.

\begin{figure}[t]
   \centering
   \includegraphics[height=2.1in, width=2.8in]{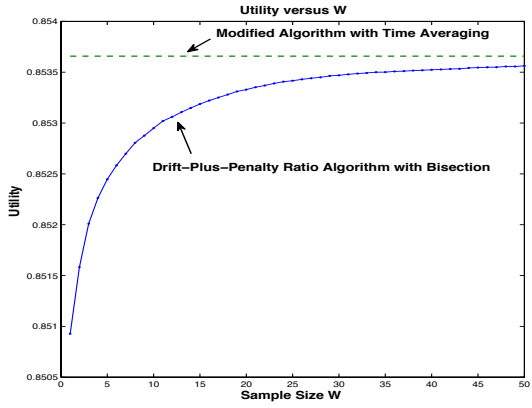} % requires the graphicx package
   \caption{Utility for the drift-plus-penalty ratio algorithm (with bisection) and the time-averaged alternative.$\vspace{-.3in}$}
   \label{fig:multiple-W}
\end{figure}

We simulate the drift-plus-penalty ratio algorithm for $10^6$ frames, 
using the bisection method with $W$ past samples of $\bv{\eta}[r]$ as
in (\ref{eq:val-theta}) of Section \ref{section:initial-info}.  We use $P^{tran} = 1.0$, $I^{max} = 5.0$. 
The vectors $\{\bv{\eta}[r]\}_{r=0}^{\infty}$ are assumed to be 
i.i.d. with independently chosen components, where $T_l^{tran}[r]$ is uniformly distributed in $[0.5, 2.5]$ for all $l$, 
and $qual_l[r]$ is uniformly distributed in $[0, l]$ for $l \in \{1, 2, 3, 4, 5\}$ (so that device $5$ tends to have the highest quality, 
while device $1$ tends to have the lowest).  We initialize $\theta_{min} = -5V$, $\theta_{max} \defequiv \sum_{l=1}^5Z_l[r]3$. 
Each step of the bisection computes $val(\theta)$ in (\ref{eq:val-theta}) 
according to a simple deterministic optimization. In particular, 
for the $w$th term in $val(\theta)$, corresponding to sample $\bv{\eta}_w$, we 
choose $Idle[r, w] = 0$ whenever $\theta \leq 0$, and $Idle[r,w] = I^{max}$ if $\theta>0$, and choose
 $l[r,w]$ as the index $l \in \{1, \ldots, 5\}$ that minimizes: 
 \[ -Vqual_{l}[r,w] + (Z_l[r]P^{tran} - \theta)T_l^{tran}[r,w]  \]

The bisection routine
is run for each frame until $\theta_{max} - \theta_{min} < 0.001$.  Using $V=100$, 
the resulting q.o.i per unit time is plotted in 
Fig. \ref{fig:multiple-W}. This increases to its optimal value as $W$ is increased.  However, in this example, 
$W$ does not need to be 
very large for accurate results:  Even $W=1$ produces a value that is near optimal (note that the $y$-axis in Fig. \ref{fig:multiple-W}
distinguishes utility only in the 3rd significant digit).

  All average power constraints are met in all simulations (for each $W$).  Results for $W=10$
are:  $q.o.i/\overline{T}  = 0.852950$, 
$\overline{T} = 3.180275$, $\overline{Idle} = 1.421260$, $\overline{y}_0 = -2.712615$, and: 
\begin{eqnarray*}
 &\overline{y}_1/\overline{T} =  0.182335 \leq 0.25& \\
&\overline{y}_2/\overline{T} =  0.249547\leq 0.25 \: \: , \: \: 
 \overline{y}_3/\overline{T} =  0.250018 \leq 0.25& \\
&\overline{y}_4/\overline{T} =  0.250032 \leq 0.25 \: \: , \: \: 
\overline{y}_5/\overline{T} =  0.250046 \leq 0.25& 
\end{eqnarray*}
It can be seen that devices $\{2, \ldots, 5\}$ are utilized to their maximum power constraints because these tend to 
give the highest quality, while average power for device 1 is slack.  

The alternative algorithm of Section \ref{section:alternative}, which does not require
a bisection routine and amounts to a simple deterministic optimization for (\ref{eq:alt-metric}) every frame,
achieves similar time average power expenditures to the above. It also achieves
utility as shown in Fig. \ref{fig:multiple-W}, being the constant that does not depend on $W$ (as
no sampling from the past is needed).  Its utility is slightly larger than that of the 
bisection algorithm, and is approached by the bisection algorithm as $W$ increases.
It appears that this algorithm is simpler and 
yields ``automatic learning''  by using the time average value $\theta[r]$, but it might have trouble adapting if
system parameters change. 

Details on the particular decisions made in the simulation on each frame are provided in Appendix G.  There, it is also shown that
the particular structure of this example admits \emph{deterministic} bounds on the constraint violations.  In particular, if 
$I^{max}$ is chosen to be suitably large (11.0 in this case), then we can guarantee that for all integers $R>0$: 
\[ \frac{\sum_{r=0}^{R-1}y_l[r]}{\sum_{r=0}^{R-1}T[r]} \leq c_l + \frac{d_1 + d_2V}{R} \]
for some positive constants $d_1, d_2$.   If we re-run the simulations using $I^{max} =11.0$ and $W=10$, 
we get similar average values 
as given above, and in particular the algorithm results in the same $\overline{Idle} \approx 1.42$, as expected. 

\section{Conclusion} 

We have developed a method for optimizing time averages in general 
renewal systems.  Every renewal frame, 
a policy is chosen that affects the frame size and also affects a penalty vector and/or an attribute vector.  
A dynamic algorithm
was developed to minimize the time average of one penalty subject to time average constraints on the others. 
A related algorithm was developed to maximize a concave function of the time average attribute vector, subject
to time average constraints. 
This work extends the theory of Lyapunov optimization to treat much more general classes of systems, including
task processing networks 
with variable length scheduling modes.

\section*{Appendix A --- Proof of Lemma \ref{lem:optimality-iid}} 

We present a definition and a simple lemma before the proof of Lemma \ref{lem:optimality-iid}. 
Define $\Gamma$ as the set of all vectors $[y_0, y_1, \ldots, y_L,  T] \in \mathbb{R}^{L+2}$ that are achievable
as time averages under i.i.d. algorithms.  Thus, a vector $[y_0, y_1, \ldots, y_L,  T]$ is in the set $\Gamma$ if and only
if there is an i.i.d. algorithm $\pi^*[r]$ such that: 
\begin{eqnarray*}
 \expect{\hat{y}_l(\pi^*[r])} &=& y_l \: \: \forall l \in \{0, 1, \ldots, L\} \\
\expect{\hat{T}(\pi^*[r])} &=& T
\end{eqnarray*}
It is easy to show that the set $\Gamma$ is bounded and convex. Further, 
for any frame $r$ and under any (possibly non-i.i.d.) algorithm $\pi[r]$, we have: 
\[ \expect{[\hat{y}_0(\pi[r]), \ldots, \hat{y}_L(\pi[r]), \hat{T}(\pi[r])]} \in \Gamma \]
This is because the policy chooses $\pi[r]$ according to some conditional distribution 
given the past history, but this distribution can be viewed on frame $r$ as one that is from an i.i.d. 
algorithm. 

\begin{lem} \label{lem:convex-hull} For any algorithm that chooses $\pi[r] \in \script{P}$ for all frames $r$, 
we have for all integers $R>0$: 
\[ \frac{1}{R}\sum_{r=0}^{R-1} \expect{[\hat{y}_0(\pi[r]), \ldots, \hat{y}_L(\pi[r]), \hat{T}(\pi[r])]} \in \Gamma \]
\end{lem} 
\begin{proof} (Lemma \ref{lem:convex-hull})  Each of the individual terms in the sum is in $\Gamma$, and so the average
of these terms is in $\Gamma$ (because $\Gamma$ is convex). 
\end{proof} 

We now prove Lemma \ref{lem:optimality-iid}.  Suppose the constraints of problem (\ref{eq:ls1})-(\ref{eq:ls3}) are feasible,
and define $ratio^{opt}$ as the infimum of the objective function over all feasible algorithms.  Fix any $\delta>0$.  Then there
must be an algorithm $\pi[r]$ that satisfies the constraints (\ref{eq:ls1})-(\ref{eq:ls3}) and yields a $\limsup$ ratio within $\delta/2$ of
the infimum.  That is, $\pi[r] \in \script{P}$ for all frames $r$, and: 
\begin{eqnarray*}
\limsup_{R\rightarrow\infty} \left[\frac{\frac{1}{R}\sum_{r=0}^{R-1}\expect{\hat{y}_0(\pi[r])}}{\frac{1}{R}\sum_{r=0}^{R-1}\expect{\hat{T}(\pi[r])}} \right] &\leq& ratio^{opt} + \delta/2 \\
\limsup_{R\rightarrow\infty} \left[\frac{\frac{1}{R}\sum_{r=0}^{R-1}\expect{\hat{y}_l(\pi[r])}}{\frac{1}{R}\sum_{r=0}^{R-1}\expect{\hat{T}(\pi[r])}}\right] &\leq& c_l \: , \: \: \forall l \in \{1, \ldots, L\} 
\end{eqnarray*}
It follows that there is a finite integer $R^*$ such that: 
\begin{eqnarray*}
\left[\frac{\frac{1}{R^*}\sum_{r=0}^{R^*-1}\expect{\hat{y}_0(\pi[r])}}{\frac{1}{R^*}\sum_{r=0}^{R^*-1}\expect{\hat{T}(\pi[r])}} \right] &\leq& ratio^{opt} + \delta  \\
 \left[\frac{\frac{1}{R^*}\sum_{r=0}^{R^*-1}\expect{\hat{y}_l(\pi[r])}}{\frac{1}{R^*}\sum_{r=0}^{R^*-1}\expect{\hat{T}(\pi[r])}}\right] &\leq& c_l + \delta \: , \: \: \forall l \in \{1, \ldots, L\} 
\end{eqnarray*}
By Lemma \ref{lem:convex-hull}, we know there is an i.i.d. algorithm $\pi^*[r]$ such that:
\begin{eqnarray*}
\frac{1}{R^*}\sum_{r=0}^{R^*-1} \expect{[\hat{y}_0(\pi[r]), \ldots, \hat{y}_L(\pi[r]), \hat{T}(\pi[r])]}  = \\
\expect{[\hat{y}_0(\pi^*[r]), \ldots, \hat{y}_L(\pi^*[r]), \hat{T}(\pi^*[r])]} 
\end{eqnarray*}
Plugging this identity into the above inequalities yields: 
\begin{eqnarray*}
\frac{\expect{\hat{y}_0(\pi^*[r])}}{\expect{\hat{T}(\pi^*[r])}} &\leq& ratio^{opt} + \delta  \\
\frac{\expect{\hat{y}_l(\pi^*[r])}}{\expect{\hat{T}(\pi^*[r])}}&\leq& c_l + \delta \: , \: \: \forall l \in \{1, \ldots, L\}  
\end{eqnarray*}
Multiplying the above by the positive number
$\expect{\hat{T}(\pi^*[r])}$ proves Lemma \ref{lem:optimality-iid}.

\section*{Appendix B --- Bound on $\overline{y}_l[R]/\overline{T}[R]$ in Theorem \ref{thm:alg1}} 

This section provides an upper bound on $\overline{y}_l[R]/\overline{T}[R]$, which 
shows how long is required to come close to meeting the 
constraints $\overline{y}_l/\overline{T} \leq c_l$.   Suppose the assumptions of Theorem \ref{thm:alg1} hold.  
From (\ref{eq:ratio-good2proof}) we have for all frames $r$: 
\begin{eqnarray*}
\Delta(\bv{Z}[r])  \leq B +T^{max}[C + Vratio^{opt} ] - Vy_0^{min}
\end{eqnarray*}
where we have used the fact that the conditional expectation of $\hat{y}_0(\pi[r])$ is bounded below by $y_0^{min}$, 
and the conditional expectation of $\hat{T}(\pi[r])$ is bounded above by $T^{max}$.  It follows that: 
\begin{equation} \label{eq:quad-drift} 
 \Delta(\bv{Z}[r]) \leq (F_1 + VF_2)/2 
 \end{equation} 
where the constants $F_1$ and $F_2$ are defined: 
\begin{eqnarray*}
F_1 \defequiv 2(B + T^{max}C) \: \: , \: \: 
F_2 \defequiv 2(T^{max} ratio^{opt} - y_0^{min})
\end{eqnarray*}
Substituting the definition of $\Delta(\bv{Z}[r])$ in (\ref{eq:quad-drift}) yields: 
\[ \expect{L(\bv{Z}[r+1]) - L(\bv{Z}[r])|\bv{Z}[r]}  \leq (F_1 + VF_2)/2 \]
Taking expectations of the above and using the law of iterated expectations yields: 
\[ \expect{L(\bv{Z}[r+1])} - \expect{L(\bv{Z}[r])} \leq (F_1 + VF_2)/2 \]
The above holds for all frames $r$.  Summing over $r \in \{0, 1, \ldots, R-1\}$ (for some positive
integer $R$) and dividing by $R$ gives: 
\[ \frac{\expect{L(\bv{Z}[R])} - \expect{L(\bv{Z}[0])}}{R} \leq (F_1 + VF_2)/2 \]
Using the fact that $L(\bv{Z}[r]) \defequiv \frac{1}{2}\sum_{l=1}^LZ_l[r]^2$,  we have: 
\[ \sum_{l=1}^L\frac{\expect{Z_l[R]^2}}{R} \leq (F_1 + VF_2) + \frac{\sum_{l=1}^L\expect{Z_l[0]^2}}{R} \]
Thus, for every $l \in \{1, \ldots, L\}$ we have: 
\[ \frac{\expect{Z_l[R]^2}}{R^2} \leq \frac{F_1 + VF_2}{R} + \frac{\sum_{l=1}^L\expect{Z_l[0]^2}}{R^2} \]
By Jensen's inequality, $\expect{Z_l[R]}^2 \leq \expect{Z_l[R]^2}$, and so for all integers
$R>0$ and all $l \in \{1, \ldots, L\}$ we have:
\[  \frac{\expect{Z_l[R]}}{R} \leq \sqrt{\frac{F_1 + VF_2}{R} + \frac{\sum_{l=1}^L\expect{Z_l[0]^2}}{R^2}} \]
Therefore, from (\ref{eq:thm1-see}) we have for all $l \in \{1, \ldots, L\}$:  
\begin{eqnarray*}
 \frac{\overline{y}_l[R]}{\overline{T}[R]} &\leq& c_l + \frac{1}{T_{min}}\frac{\expect{Z_l[R]}}{R} \\
 &\leq& 
 c_l + \frac{1}{T_{min}} \sqrt{\frac{F_1 + VF_2}{R} + \frac{\sum_{l=1}^L\expect{Z_l[0]^2}}{R^2}} 
 \end{eqnarray*}

\section*{Appendix C --- Analysis of Alternative Formulation} 

Consider the alternative problem from Section \ref{section:alt-form}: 
\begin{eqnarray}
\mbox{Minimize:} & \limsup_{R\rightarrow\infty} \overline{y}_0[R] \label{eq:alt1} \\
\mbox{Subject to:} & \limsup_{R\rightarrow\infty}\frac{\overline{y}_l[R]}{\overline{T}[R]}\leq c_l \: \: \forall l \in \{1, \ldots, L\} \label{eq:alt2}  \\
& \pi[r] \in \script{P} \: \: \forall r \in \{0, 1,2, \ldots \} \label{eq:alt3} 
\end{eqnarray}
Assume the same boundedness assumptions of Section \ref{section:boundedness} hold. Further
assume that a $C$-additive approximation of
the algorithm of Section \ref{section:alt-form} is implemented, 
so that every frame $r$ we observe $\bv{Z}[r]$ and choose $\pi[r] \in \script{P}$
to yield: 
\begin{eqnarray} 
\expect{V\hat{y}_0(\pi[r]) + \sum_{l=1}^LZ_l[r][\hat{y}_l(\pi[r]) - c_l\hat{T}(\pi[r])]|\bv{Z}[r]} \leq \nonumber \\
C + \expect{V\hat{y}_0(\pi^*[r]) + \sum_{l=1}^LZ_l[r][\hat{y}_l(\pi^*[r]) - c_l\hat{T}(\pi^*[r])]} \label{eq:capp-f1}
\end{eqnarray} 
where $C$ is a given non-negative constant, and $\pi^*[r]$ is any i.i.d. algorithm. 
\begin{thm} Under the above assumptions, and assuming the constraints of the 
problem (\ref{eq:alt1})-(\ref{eq:alt3}) are feasible, then: 

a) For all $l \in \{1, \ldots, L\}$ we have: 
\begin{eqnarray}
 \limsup_{R\rightarrow\infty} \overline{y}_l[R]/\overline{T}[R] &\leq& c_l \: \: \forall l \in \{1, \ldots, L\} \\
 \limsup_{R\rightarrow\infty} \frac{\sum_{r=0}^{R-1} y_l[r]}{\sum_{r=0}^{R-1}T[r]} &\leq& c_l \: \: \: (w.p.1)  
\end{eqnarray}

b) For all integers $R>0$ we have: 
\begin{eqnarray*}
\overline{y}_0[R] \leq y_0^{opt} + \frac{B + C}{V} + \frac{\expect{L(\bv{Z}[0])}}{VR}
\end{eqnarray*}
where $B$ is defined in (\ref{eq:B}), and $y_0^{opt}$ is the infimum value of (\ref{eq:alt1}) subject to (\ref{eq:alt2})-(\ref{eq:alt3}). 
\end{thm} 
 
\begin{proof} 
Consider any frame $r \in \{0, 1, 2, \ldots\}$. 
From (\ref{eq:drift}) we have: 
\begin{eqnarray*}
&\Delta(\bv{Z}[r]) + V\expect{\hat{y}_0(\pi[r])|\bv{Z}[r]} \leq B \\
&+ \expect{V\hat{y}_0(\pi[r]) + \sum_{l=1}^LZ_l[r][\hat{y}_l(\pi[r]) - c_l\hat{T}(\pi[r])]|\bv{Z}[r]} 
\end{eqnarray*}
Substituting (\ref{eq:capp-f1}) yields:
\begin{eqnarray}
\Delta(\bv{Z}[r]) + V\expect{\hat{y}_0(\pi[r])|\bv{Z}[r]} \leq B + C \nonumber  \\
+ \expect{V\hat{y}_0(\pi^*[r]) + \sum_{l=1}^LZ_l[r][\hat{y}_l(\pi^*[r]) - c_l\hat{T}(\pi^*[r])]}  \label{eq:afp1} 
\end{eqnarray}
where $\pi^*[r]$ is any i.i.d. algorithm. 
As in Lemma \ref{lem:optimality-iid}, it can be shown that if the problem (\ref{eq:alt1})-(\ref{eq:alt3}) is feasible, then 
for any $\delta>0$ there is an i.i.d. algorithm $\pi^*[r]$ such that: 
\begin{eqnarray*}
\expect{\hat{y}_0(\pi^*[r])} &\leq& y_0^{opt} + \delta \\
\expect{\hat{y}_l(\pi^*[r])}/\expect{\hat{T}(\pi^*[r])} &\leq& c_l + \delta \: \: \forall l \in \{1, \ldots, L\}
\end{eqnarray*}
Plugging the above into the right-hand-side of (\ref{eq:afp1})  yields: 
\begin{eqnarray*}
\Delta(\bv{Z}[r]) + V\expect{\hat{y}_0(\pi[r])|\bv{Z}[r]} \leq B + C + V(y_0^{opt}+\delta) \nonumber \\
+ \sum_{l=1}^LZ_l[r]\left[(c_l+\delta)\expect{\hat{T}(\pi^*[r])}  - \expect{c_l\hat{T}(\pi^*[r])}\right]
\end{eqnarray*}
Taking $\delta\rightarrow 0$ yields: 
\begin{eqnarray}
\Delta(\bv{Z}[r]) + V\expect{\hat{y}_0(\pi[r])|\bv{Z}[r]} \leq B + C + Vy_0^{opt} \label{eq:afp3} 
\end{eqnarray}
From this we obtain: 
\[ \Delta(\bv{Z}[r]) \leq B + C + V(y_0^{opt} - y_0^{min}) \]
Thus, the quadratic Lyapunov drift is less than or equal to a constant, from which we obtain the result of part (a) by the
same argument as in the proof of Theorem \ref{thm:alg1}. 

To prove part (b), taking expectations of (\ref{eq:afp3}) yields: 
\[ \expect{L(\bv{Z}[r+1])} - \expect{L(\bv{Z}[r])} + V\expect{y_0[r]} \leq B + C + Vy_0^{opt} \]
Fix an integer $R>0$. Summing the above 
over $r \in \{0, \ldots R-1\}$ and dividing by $R$ gives: 
\[ \frac{\expect{L(\bv{Z}[R])} - \expect{L(\bv{Z}[0])}}{R} + V\overline{y}_0[R] \leq B + C + Vy_0^{opt} \]
Rearranging terms and noting that $\expect{L(\bv{Z}[R])}\geq 0$ proves the result. 
\end{proof} 

\section*{Appendix D --- Analysis of the Alternative Algorithm with Time Averaging} 

Here we consider the original problem of 
minimizing $\overline{y}_0/\overline{T}$ subject
to $\overline{y}_l/\overline{T} \leq c_l$ for $l \in\{1, \ldots, L\}$, and
analyze the alternative algorithm (with time averaging) described in 
Section \ref{section:alternative}.  Recall that queues $\bv{Z}[r]$ still operate according
to (\ref{eq:z-update}). 
Define $\theta[0] \defequiv 0$, and for integers $r>0$ define
$y_0^{av}[r]$, $T^{av}[r]$, $\theta[r]$ by: 
\begin{eqnarray*}
y_0^{av}[r] \defequiv \frac{1}{r}\sum_{i=0}^{r-1} y_0[r] \: \: , \: \: 
T^{av}[r] \defequiv \frac{1}{r}\sum_{i=0}^{r-1} T[r] \: \: , \: \: 
\theta[r] \defequiv \frac{y_0^{av}[r]}{T^{av}[r]}
\end{eqnarray*}
Assume that we use a $C$-additive approximation for the policy selection step in (\ref{eq:alt-metric}), so that
every frame $r$ we observe $\bv{Z}[r]$ and $\theta[r]$ and choose $\pi[r] \in \script{P}$ to yield: 
\begin{eqnarray}
\mathbb{E}\{V[\hat{y}_0(\pi[r]) - \theta[r]\hat{T}(\pi[r])]|\bv{Z}[r], \theta[r]\} \nonumber \\
  +  \mathbb{E}\{\sum_{l=1}^LZ_l[r][\hat{y}_l(\pi[r]) - c_l\hat{T}(\pi[r])]|\bv{Z}[r], \theta[r]\} \leq \nonumber \\
 C + V\expect{\hat{y}_0(\pi^*[r])} - \theta[r]\expect{\hat{T}(\pi^*[r])} \nonumber \\
 + \sum_{l=1}^LZ_l[r][\expect{\hat{y}_l(\pi^*[r])} - c_l\expect{\hat{T}(\pi^*[r])}] \label{eq:dpp-alt1} 
\end{eqnarray}

We now make the following convergence assumption: 

\emph{Assumption 1:} There are constants $y_0^{av}$, $T^{av}$, and $\theta^* \defequiv y_0^{av}/T^{av}$, such that 
under the implementation we have the following convergence properties: 
\begin{eqnarray}
\lim_{r\rightarrow\infty} (y_0^{av}[r], T^{av}[r], \theta[r]) = (y_0^{av}, T^{av}, \theta^*) \: \: (w.p.1) \label{eq:ass1-a} \\
\lim_{r\rightarrow\infty} (\expect{y_0^{av}[r]}, \expect{T^{av}[r]}, \expect{\theta[r]}) = (y_0^{av}, T^{av}, \theta^*) \label{eq:ass1-b} 
\end{eqnarray}

The equality (\ref{eq:ass1-b}) typically holds whenever (\ref{eq:ass1-a}) holds.  
Indeed, taking an expectation of (\ref{eq:ass1-a}) and assuming we can pass the expectation
through the limit yields (\ref{eq:ass1-b}).  We can exchange the limit and the expectation 
whenever $y_0^{av}[r]$, $T^{av}[r]$, $\theta[r]$ are deterministically bounded  for all $r$ and all sample paths, 
or when milder conditions
hold that allow the Lebesgue dominated convergence theorem to be applied \cite{williams-martingale}. 

If (\ref{eq:ass1-a}) holds, it is easy to show that, with probability 1: 
\[ \lim_{R\rightarrow\infty} \frac{1}{R}\sum_{r=0}^{R-1}\theta[r]T[r] = \theta^*T^{av} = y_0^{av} \]
It is natural (and useful) to assume the above also holds in expectation: 

\emph{Assumption 2:} With $\theta^*$, $T^{av}$, $y_0^{av}$ as defined in Assumption 1, we have:
\[ \lim_{R\rightarrow\infty} \frac{1}{R}\sum_{r=0}^{R-1}\expect{\theta[r]T[r]} = \theta^*T^{av} = y_0^{av} \]

\begin{thm} \label{thm:alt-time-av} Assume the constraints of the problem (\ref{eq:ls1})-(\ref{eq:ls3}) are feasible, 
that $\expect{L(\bv{Z}[0])} < \infty$, and that we use
a $C$-additive approximation as described above every frame. Then: 

a) For all $l \in \{1, \ldots, L\}$ we have: 
\begin{eqnarray*}
 \limsup_{R\rightarrow\infty} \overline{y}_l[R]/\overline{T}[R] &\leq& c_l \: \: \forall l \in \{1, \ldots, L\} \\
 \limsup_{R\rightarrow\infty} \frac{\sum_{r=0}^{R-1} y_l[r]}{\sum_{r=0}^{R-1}T[r]} &\leq& c_l \: \: \: (w.p.1) 
\end{eqnarray*}

b) If Assumptions 1 and 2 hold, then 
the achieved value of $\overline{y}_0/\overline{T}$ satisfies the following with probability 1: 
\[ \lim_{R\rightarrow\infty} y_0^{av}[R]/T^{av}[R] \leq ratio^{opt} + (B + C)/VT^{min} \]
where $B$ is defined as in (\ref{eq:B}) (with a conditional expectation that is also given $\theta[r]$). 
\end{thm} 

\begin{proof} 
Define $\Delta(\bv{Z}[r], \theta[r])$ as the conditional drift, conditioned also on knowledge of $\theta[r]$.
This has the same form as (\ref{eq:drift}), and so from (\ref{eq:drift}) we have: 
\begin{eqnarray*}
\Delta(\bv{Z}[r], \theta[r]) + V\expect{\hat{y}_0(\pi[r]) - \theta[r]\hat{T}(\pi[r])|\bv{Z}[r], \theta[r]} \leq \\
B + V\expect{\hat{y}_0(\pi[r]) - \theta[r]\hat{T}(\pi[r])|\bv{Z}[r], \theta[r]} \\
+ \expect{\sum_{l=1}^LZ_l[r][\hat{y}_l(\pi[r]) - c_l\hat{T}(\pi[r])]|\bv{Z}[r], \theta[r]} 
\end{eqnarray*}
Using (\ref{eq:dpp-alt1}) in the right-hand-side above yields: 
\begin{eqnarray*}
\Delta(\bv{Z}[r], \theta[r]) + V\expect{\hat{y}_0(\pi[r]) - \theta[r]\hat{T}(\pi[r])|\bv{Z}[r], \theta[r]} \leq \\
B + C + V\expect{\hat{y}_0(\pi^*[r])} - \theta[r]\expect{\hat{T}(\pi^*[r])} \\
+ \sum_{l=1}^LZ_l[r][\expect{\hat{y}_l(\pi^*[r])} - c_l\expect{\hat{T}(\pi^*[r])}] 
\end{eqnarray*}
where $\pi^*[r]$ is any i.i.d. algorithm.  Now plug the i.i.d. algorithm $\pi^*[r]$ from (\ref{eq:iid1})-(\ref{eq:iid2})
into the right-hand-side of the above and take $\delta\rightarrow0$ to yield: 
\begin{eqnarray}
\Delta(\bv{Z}[r], \theta[r]) + V\expect{\hat{y}_0(\pi[r]) - \theta[r]\hat{T}(\pi[r])|\bv{Z}[r], \theta[r]} \leq \nonumber \\
B + C + Vratio^{opt}\expect{\hat{T}(\pi^*[r])} - V\theta[r]\expect{\hat{T}(\pi^*[r])}  \label{eq:alt33} 
\end{eqnarray}
The above has the form: 
\[ \Delta(\bv{Z}[r], \theta[r]) \leq F\]
for some finite constant $F$, and so part (a) holds by arguments similar to those in the proof of Theorem \ref{thm:alg1}. 

Now taking expectations of (\ref{eq:alt33}) yields: 
\begin{eqnarray*}
\expect{L(\bv{Z}[r+1])} - \expect{L(\bv{Z}[r])}  \nonumber \\
+ V[\expect{\hat{y}_0[r]} - \expect{\theta[r]T[r]}] \leq \nonumber \\
B + C + Vratio^{opt}\expect{\hat{T}(\pi^*[r])} \nonumber \\
- V\expect{\theta[r]}\expect{\hat{T}(\pi^*[r])} 
\end{eqnarray*}
Summing the above over $r \in \{0, \ldots, R-1\}$ and dividing by $R$ yields: 
\begin{eqnarray*}
 \frac{\expect{L(\bv{Z}[R])} - \expect{L(\bv{Z}[0])}}{R} \\
 + V[\overline{y}_0[R] - \frac{1}{R}\sum_{r=0}^{R-1}\expect{\theta[r]T[r]}] \leq \\
 B + C + Vratio^{opt}\expect{\hat{T}(\pi^*[r])} \\
 - V\expect{\hat{T}(\pi^*[r])}\frac{1}{R}\sum_{r=0}^{R-1}\expect{\theta[r]}
\end{eqnarray*}
Thus: 
\begin{eqnarray*}
 V[\overline{y}_0[R] - \frac{1}{R}\sum_{r=0}^{R-1}\expect{\theta[r]T[r]}] \leq \\
 B + C + Vratio^{opt}\expect{\hat{T}(\pi^*[r])} \\
 - V\expect{\hat{T}(\pi^*[r])}\frac{1}{R}\sum_{r=0}^{R-1}\expect{\theta[r]} + \frac{\expect{L(\bv{Z}[0])}}{R}
\end{eqnarray*}
However, by Assumptions 1 and 2: 
\[ \lim_{R\rightarrow\infty} [\overline{y}_0[R] - \frac{1}{R}\sum_{r=0}^{R-1}\expect{\theta[r]T[r]}] = y_0^{av} - y_0^{av} = 0 \]
Thus, taking limits of the above yields: 
\begin{eqnarray*}
0\leq B + C + Vratio^{opt}T^* - VT^*\theta^* 
\end{eqnarray*}
where we have defined $T^*\defequiv \expect{\hat{T}(\pi^*[r])}$, and we have used the fact that 
if $\expect{\theta[r]} \rightarrow \theta^*$, then so does its time average.
Rearranging terms yields: 
\[ y_0^{av}/T^{av} \defequiv \theta^* \leq ratio^{opt} + \frac{B+C}{VT^*} \leq ratio^{opt} + \frac{B+C}{VT^{min}} \]
This proves Theorem \ref{thm:alt-time-av}.
\end{proof} 

\section*{Appendix E --- Variation on Jensen's Inequality} 

Here we prove part (a) of Lemma \ref{lem:jensen}.  Recall that $(T, \bv{\gamma})$ is a random vector
with arbitrary joint distribution such that $T>0$ and $\bv{\gamma} \in \script{R}$ with probability 1, and 
$0 < \expect{T} < \infty$. 
Define $\{(T[r],\bv{\gamma}[r])\}_{r=0}^{\infty}$ as an infinite sequence of i.i.d. random vectors, 
 each distributed the same as $(T, \bv{\gamma})$.  Define $t[0] = 0$, and for integers $r>0$ define $t[r] \defequiv \sum_{i=0}^{r-1}T[i]$. 
 The value of $t[r]$ can be viewed as the $r$th \emph{renewal time} in a process with i.i.d. renewal durations of size $T[r]$. 
 Define $\bv{\gamma}(t)$ as a random vector process defined 
 over continuous time $t \geq 0$, taking the value $\bv{\gamma}[r]$ whenever $t$ is in the $r$th renewal interval, so that: 
 \[ \bv{\gamma}(t) = \bv{\gamma}[r] \: \: \mbox{ if and only if $t[r] \leq t < t[r+1]$} \]
 Therefore, for any integer $r>0$: 
 \[ \frac{1}{t[r]}\int_{0}^{t[r]} \phi(\bv{\gamma}(t)) dt  = \frac{\sum_{i=0}^{r-1}T[i]\phi(\bv{\gamma}[i])}{\sum_{i=0}^{r-1}T[i]} \]
 On the other hand, by Jensen's inequality for integrals of  concave functions, we have: 
 \begin{eqnarray*}
 \frac{1}{t[r]}\int_{0}^{t[r]} \phi(\bv{\gamma}(t)) dt &\leq& \phi\left(\frac{1}{t[r]}\int_{0}^{t[r]}\bv{\gamma}(t)dt\right) \\
 &=& \phi\left(\frac{\sum_{i=0}^{r-1}\bv{\gamma}[i]T[i]}{\sum_{i=0}^{r-1}T[i]}  \right)
 \end{eqnarray*}
 Thus, for all integers $r>0$ we have: 
 \begin{equation} \label{eq:jensen-proof-dude} 
  \frac{\frac{1}{r}\sum_{i=0}^{r-1}T[i]\phi(\bv{\gamma}[i])}{\frac{1}{r}\sum_{i=0}^{r-1}T[i]} \leq \phi\left(\frac{\frac{1}{r}\sum_{i=0}^{r-1}\bv{\gamma}[i]T[i]}{\frac{1}{r}\sum_{i=0}^{r-1}T[i]}  \right) 
  \end{equation} 
 However, by the law of large numbers, we know: 
 \begin{eqnarray*}
 \lim_{r\rightarrow\infty} \frac{1}{r}\sum_{i=0}^{r-1}T[i] &=& \expect{T} \: \: (w.p.1) \\
  \lim_{r\rightarrow\infty} \frac{1}{r}\sum_{i=0}^{r-1}T[i]\phi(\bv{\gamma}[i]) &=& \expect{T\phi(\bv{\gamma})} \: \: (w.p.1) \\
 \lim_{r\rightarrow\infty} \frac{1}{r}\sum_{i=0}^{r-1}\bv{\gamma}[i]T[i] &=& \expect{T\bv{\gamma}} \: \: (w.p.1) 
 \end{eqnarray*}
 Taking limits of the above as $r\rightarrow\infty$ in (\ref{eq:jensen-proof-dude}) and using the above identities 
 together with continuity of $\phi(\bv{\gamma})$ proves that: 
 \[ \frac{\expect{T\phi(\bv{\gamma})}}{\expect{T}} \leq \phi\left(\frac{\expect{T\bv{\gamma}}}{\expect{T}}\right) \]
 
\section*{Appendix F --- Analysis of Utility Optimization} 

Here we prove Theorem \ref{thm:utility}. For simplicity, assume that the set $\Gamma$ of all expectations 
$[\bv{y}[r], \bv{x}[r], T[r]]$ under i.i.d. algorithms is closed.  Then, similar
 to Lemma \ref{lem:optimality-iid}, it can be shown that if the problem (\ref{eq:tran1})-(\ref{eq:tran3}) is feasible, then there is 
 an i.i.d. algorithm $\pi^*[r]$ and a vector $\bv{\gamma} \in \script{R}$ such that: 
 \begin{eqnarray}
 \phi(\bv{\gamma}^*) &=& util^{opt} \label{eq:iid-t1}  \\
 \frac{\expect{\hat{x}_m(\pi^*[r])}}{\expect{\hat{T}(\pi^*[r])}} &=& \gamma_m^* \: \: \forall m \in \{1, \ldots, M\} \label{eq:iid-t2}  \\
  \frac{\expect{\hat{y}_l(\pi^*[r])}}{\expect{\hat{T}(\pi^*[r])}} &\leq& c_l  \: \: \forall l \in \{1, \ldots, L\} \label{eq:iid-t3} 
 \end{eqnarray}
 If the set $\Gamma$ is not closed, then (\ref{eq:iid-t1})-(\ref{eq:iid-t3}) can be modified to show they hold to within any $\delta>0$, 
 and the same result we derive below can be recovered by taking 
 $\delta\rightarrow 0$,  as in the proof of Theorem \ref{thm:alg1}. 
 
Now define $\bv{Q}[r] \defequiv [\bv{Z}[r]; \bv{G}[r]]$ as the collection of all queues, and 
 define the Lyapunov function: 
 \[ L(\bv{Q}[r]) \defequiv \frac{1}{2}\sum_{l=1}^LZ_l[r]^2 + \frac{1}{2}\sum_{m=1}^MG_m[r]^2 \]
 For simplicity, assume initial conditions satisfy $Z_l[0] = G_m[0] = 0$ for all $l$ and $m$. 
 Define $\Delta(\bv{Q}[r]) \defequiv \expect{L(\bv{Q}[r+1]) - L(\bv{Q}[r])|\bv{Q}[r]}$. 
 By an argument similar to that given in Lemma \ref{lem:drift-compute}, we can square the queue update
 equations (\ref{eq:z-update}) and (\ref{eq:g-update}) to obtain the following drift-plus-penalty bound: 
 \begin{eqnarray*}
 \Delta(\bv{Q}[r]) - V\expect{T[r]\phi(\bv{\gamma}[r])|\bv{Q}[r]} \leq \\
 D - V\expect{T[r]\phi(\bv{\gamma}[r])|\bv{Q}[r]} \\
  + \expect{\sum_{l=1}^LZ_l[r][y_l[r] - c_lT[r]]|\bv{Q}[r]} \\
  + \expect{\sum_{m=1}^M G_m[r][T[r]\gamma_m[r] - x_m[r]] |\bv{Q}[r]} 
 \end{eqnarray*}
where $D$ is defined in (\ref{eq:D}). 
  This drift-plus-penalty bound can be rearranged as follows: 
   \begin{eqnarray}
 \Delta(\bv{Q}[r]) - V\expect{T[r]\phi(\bv{\gamma}[r])|\bv{Q}[r]} \leq \nonumber \\
 D - \expect{\left[V\phi(\bv{\gamma}[r]) - \sum_{m=1}^MG_m[r]\gamma_m[r]\right]T[r]|\bv{Q}[r]} \nonumber \\
  + \expect{\sum_{l=1}^LZ_l[r]y_l[r] - \sum_{m=1}^M G_m[r]x_m[r] |\bv{Q}[r]} \nonumber \\
  -\sum_{l=1}^LZ_l[r]c_l\expect{T[r]|\bv{Q}[r]}  \label{eq:util-dpp} 
 \end{eqnarray}
  However, by the auxiliary variable update algorithm, we know that for every frame $r$: 
  \[ V\phi(\bv{\gamma}[r]) - \sum_{m=1}^MG_m[r]\gamma_m[r] \geq V\phi(\bv{\gamma}^*) - \sum_{m=1}^MG_m[r]\gamma_m^* \]
  for any vector $\bv{\gamma}^* = (\gamma_1^*, \ldots, \gamma_M^*) \in \script{R}$. 
  Further, because we use a $C$-additive approximation when choosing $\pi[r] \in \script{P}$,  every frame $r$ we have:  
 \begin{eqnarray*} 
 \frac{\expect{ \sum_{l=1}^LZ_l[r]\hat{y}_l(\pi[r]) - \sum_{m=1}^MG_m[r]\hat{x}_m(\pi[r])|\bv{Q}[r]}}{\expect{\hat{T}(\pi[r])|\bv{Q}[r]}}  \leq C  \nonumber \\
+ \frac{\sum_{l=1}^LZ_l[r]\expect{\hat{y}_l(\pi^*[r])} - \sum_{m=1}^MG_m[r]\expect{\hat{x}_m(\pi^*[r])}}{\expect{\hat{T}(\pi^*[r])}} 
 \end{eqnarray*} 
 where $\pi^*[r]$ is any i.i.d. algorithm.  Plugging the above two inequalities into the right-hand-side of (\ref{eq:util-dpp}) yields: 
 \begin{eqnarray}
 \Delta(\bv{Q}[r]) - V\expect{T[r]\phi(\bv{\gamma}[r])|\bv{Q}[r]} \leq \nonumber \\
 D  - \left[V\phi(\bv{\gamma}^*) - \sum_{m=1}^MG_m[r]\gamma_m^*\right]\expect{T[r]|\bv{Q}[r]} \nonumber \\
  + C\expect{T[r]|\bv{Q}[r]} + \expect{T[r]|\bv{Q}[r]} \times \nonumber \\
  \left[\sum_{l=1}^LZ_l[r]\frac{\expect{\hat{y}_l(\pi^*[r])}}{\expect{\hat{T}(\pi^*[r])}} - \sum_{m=1}^M G_m[r]\frac{\expect{\hat{x}_m(\pi^*[r])}}{\expect{\hat{T}(\pi^*[r])}}\right] \nonumber \\
  -\sum_{l=1}^LZ_l[r]c_l\expect{T[r]|\bv{Q}[r]}  \label{eq:util-dpp2} 
 \end{eqnarray}
 where $\bv{\gamma}^*$ is any vector in $\script{R}$, and $\pi^*[r]$ is any i.i.d. algorithm. 
 Plugging the vector $\bv{\gamma}^*$ and the i.i.d. algorithm $\pi^*[r]$ from (\ref{eq:iid-t1})-(\ref{eq:iid-t3}) into the right-hand-side
 of (\ref{eq:util-dpp2}) yields: 
 \begin{eqnarray}
 \Delta(\bv{Q}[r]) - V\expect{T[r]\phi(\bv{\gamma}[r])|\bv{Q}[r]} \leq \nonumber \\
 D + C\expect{T[r]|\bv{Q}[r]}  - Vutil^{opt}\expect{T[r]|\bv{Q}[r]}  \label{eq:util3} 
 \end{eqnarray}
 It follows from the above that $\Delta(\bv{\Theta}[r]) \leq F$ 
 for some constant $F$.  Thus, we know all queues $Z_l[r]$ and $G_m[r]$ are mean rate stable and rate stable \cite{lyap-opt2}. 
This proves  part (a) by an argument similar to that given in the proof of Theorem \ref{thm:alg1}. 

To prove part (b), we have by taking expectations of (\ref{eq:util3}): 
\begin{eqnarray*}
 \expect{L(\bv{Q}[r+1])} - \expect{L(\bv{Q}[r])}  - V\expect{T[r]\phi(\bv{\gamma}[r])} \leq \nonumber \\
 D + C\expect{T[r]}  - Vutil^{opt}\expect{T[r]}  
 \end{eqnarray*}
 Summing the above over $r \in \{0, \ldots, R-1\}$, dividing by $RV$, recalling our assumption that $L(\bv{Q}[0]) = 0$, 
 removing the non-negative
 term $\expect{L(\bv{Q}[R])}$, and rearranging terms yields: 
 \begin{eqnarray*}
 \frac{1}{R}\sum_{r=0}^{R-1}\expect{T[r]\phi(\bv{\gamma}[r])} \geq 
  util^{opt}\overline{T}[r] - \frac{D+C\overline{T}[R]}{V}
 \end{eqnarray*}
 where we have used $\overline{T}[R] \defequiv \frac{1}{R}\sum_{r=0}^{R-1}\expect{T[r]}$.  Dividing
 both sides by $\overline{T}[R]$ and using the variation on Jensen's inequality
 (Lemma \ref{lem:jensen}) yields: 
 \begin{equation} \label{eq:use-util}
  \phi\left(\frac{\frac{1}{R}\sum_{r=0}^{R-1}\expect{T[r]\bv{\gamma}[r]}}{\overline{T}[R]}\right) \geq   util^{opt}  - \frac{D}{V\overline{T}[R]} - \frac{C}{V}
  \end{equation} 
The above holds for all integers $R>0$. 

However, we have from (\ref{eq:g-update}) (similar to derivation
of (\ref{eq:thm1-see})): 
\[ \frac{\expect{\bv{G}[R]}}{R} + \overline{\bv{x}}[R] \geq \frac{1}{R}\sum_{r=0}^{R-1}\expect{T[r]\bv{\gamma}[r]} \]
where the above vector inequality is taken entrywise.   Using this in (\ref{eq:use-util}) together with the fact that 
$\phi(\bv{\gamma})$ is entrywise non-decreasing, we have: 
 \begin{equation*}
  \phi\left(\frac{\overline{\bv{x}}[R] + \expect{\bv{G}[R]}/R}{\overline{T}[R]}\right) \geq util^{opt} -  \frac{D}{V\overline{T}[R]} - \frac{C}{V}
  \end{equation*} 
  Using the fact that $\overline{T}[R] \geq T^{min}$ gives: 
  \begin{equation} \label{eq:use-util2}
  \phi\left(\frac{\overline{\bv{x}}[R] + \expect{\bv{G}[R]}/R}{\overline{T}[R]}\right) \geq util^{opt}  -  \frac{D}{VT^{min}} - \frac{C}{V}
  \end{equation} 

Now recall that all queues $G_m[r]$ are mean rate
stable, so that: 
\[ \lim_{R\rightarrow\infty} \frac{\expect{\bv{G}[R]}}{R} = \bv{0} \]
Taking a $\liminf$ of (\ref{eq:use-util2}) as $R\rightarrow\infty$
and using continuity of $\phi(\cdot)$ proves part (b) of Theorem \ref{thm:utility}.

\section*{Appendix G --- Simulation Details and Deterministic Bounds} 

Here we provide simulation details for the particular system of Section \ref{section:simulation}.  We also show
that deterministic bounds on the constraint violations are computable if the $I^{max}$ parameter is chosen suitably
large. 

\subsection{Simulation Details} 

The bisection algorithm for the above example was implemented by computing $val(\theta)$ in 
(\ref{eq:val-theta}), where the term for each sample $\bv{\eta}_w$ is found by observing $\bv{Z}[r]$
and choosing $l[r] \in \{1, \ldots, 5\}$
and $Idle[r] \in [0, I^{max}]$ to minimize the following deterministic expression: 
\begin{eqnarray*}
 -Vqual_{l[r]}[r] + \sum_{l=1}^LZ_l[r][0.5 + P^{tran}T_l^{tran}[r]1_{\{l[r]=l\}}]  \\
 - \theta(0.5 + T_{l[r]}^{tran}[r] + Idle[r])
 \end{eqnarray*}
 This is solved by choosing $Idle[r] = 0$ whenever $\theta \leq 0$, and $Idle[r] = I^{max}$ if $\theta>0$, and choosing
 $l[r]$ as the index $l \in \{1, \ldots, 5\}$ that minimizes: 
 \[ -Vqual_{l}[r] + (Z_l[r]P^{tran} - \theta)T_l^{tran}[r]  \]
 
 The alternative algorithm with time averaging is implemented by observing $\bv{Z}[r]$ and $\theta[r]$ and 
 minimizing 
 (\ref{eq:alt-metric}), which in this context amounts to choosing $l[r] \in \{1, \ldots, 5\}$, $Idle[r] \in [0, I^{max}]$ to minimize: 
 \begin{eqnarray*}
 -Vqual_{l[r]}[r] - V\theta[r](T_{l[r]}^{tran}[r] + Idle[r]) +\\
  \sum_{l=1}^LZ_l[r][P^{tran}T_l^{tran}[r]1_{\{l[r]=l\}} - c_l(T_{l[r]}^{tran}[r] + Idle[r])]
 \end{eqnarray*}
 This amounts to choosing $Idle[r] = 0$ whenever $V\theta[r] + \sum_{l=1}^LZ_l[r]c_l \leq 0$, and $Idle[r] = I^{max}$ else, 
 and choosing $l[r]$ as the index $l \in \{1, \ldots, 5\}$ that minimizes: 
 \[ -Vqual_l[r] - T_l^{tran}[r]\left[V\theta[r] - Z_l[r]P^{tran} +  \sum_{k=1}^LZ_k[r]c_k\right] \]

\subsection{Deterministic Queue Bounds} \label{section:structure}

Here we show that, if $I^{max}$ is chosen to be suitably large, then 
the drift-plus-penalty ratio algorithm for this context yields 
\emph{deterministic} bounds on $Z_l[r]$.  The ratio to minimize is: 
\[ \frac{\expect{-Vqual_{l[r]}[r] + \sum_{l=1}^LZ_l[r]\hat{y}_l(\pi[r])  |\bv{Z}[r]}}{\expect{0.5 + T_{l[r]}^{tran}[r] + Idle[r] |\bv{Z}[r]}} \]
Because $\hat{y}_l(\pi[r]) \geq 0.5$ for all $l$ and all policy choices, 
and $-Vqual_{l[r]}[r] \geq -5V$ for all policy choices,  the numerator above is positive whenever: 
\[ \sum_{l=1}^L Z_l[r](0.5) > -5V \]
In particular, the algorithm chooses $Idle[r] = I^{max}$ whenever $Z_l[r] > 10V$ for any queue $l \in \{1, \ldots, L\}$. 

Recall that the $Z_l[r]$ update is given by: 
\[ Z_l[r+1] = \max[Z_l[r] +\hat{y}_l(\pi[r]) - 0.25\hat{T}(\pi[r]), 0] \] 
Because $\hat{y}_l(\pi[r]) \leq 0.5 + 2.5P^{tran}$, and $\hat{T}(\pi[r])  \geq 1.0 + Idle[r]$, whenever $Idle[r] = I^{max}$ we have: 
\[ Z_l[r+1] \leq \max[Z_l[r] + 0.5  + 2.5P^{tran}- 0.25(1.0 + I^{max}), 0] \]
Therefore, $Z_l[r]$ cannot increase on the next slot if $Idle[r] = I^{max}$ and if: 
\[ I^{max} \geq \frac{0.5 + 2.5P^{tran}- 0.25}{0.25} = 1.0 + 10P^{tran} \]
This means that if any queue $Z_l[r]$ exceeds $10V$, then $Idle[r] = I^{max}$ and 
it cannot increase further. Because the maximum increase in $Z_l[r]$ is $0.5 + 2.5P^{tran} - 0.25$,  
for all $l\in \{1,\ldots, 5\}$ and 
all frames $r$ we have: 
\[ 0 \leq Z_l[r] \leq 10V + 0.5 + 2.5P^{tran} - 0.25 = d_1V + d_2\]
provided that this inequality holds for $r=0$, where $d_1 \defequiv 10$ and $d_2 \defequiv 0.25+ 2.5P^{tran}$.   
Thus, from (\ref{eq:det-bound}) we have for all $l \in \{1, \ldots, L\}$
and all integers $R>0$: 
\begin{eqnarray}
 \frac{\sum_{r=0}^{R-1}y_l[r]}{\sum_{r=0}^{R-1}T[r]} &\leq& c_l + \frac{d_1V + d_2}{\sum_{r=0}^{R-1}T[r]} \nonumber \\
 &\leq& c_l +  \frac{d_1V + d_2}{R} \label{eq:see-structure} 
\end{eqnarray}
where we have used the fact that $T[r] \geq 1.0$ for all $r$.  This provides a deterministic guarantee on the worst-case constraint
violation over any interval of frames starting at frame $0$.  The deviation is $O(1/R)$, which decays 
at a faster rate than the general $O(\sqrt{1/R})$ bound in (\ref{eq:obtained-bound}). 

\section*{Appendix H --- Convergence Under a Slater Condition} 

We first present a definition and two theorems from \cite{lyap-opt2}. Define $\script{H}[0]\defequiv\bv{Z}[0]$, and for integers
$r>0$ define $\script{H}[r]$ as the system history up to frame $r$, being the queue values up to and including
frame $r$, and the 
penalties up to but not including frame $r$: 
\[ \script{H}[r] \defequiv [\bv{Z}[0], \bv{Z}[1], \ldots, \bv{Z}[r], \bv{y}[0], \bv{y}[1], \ldots, \bv{y}[r-1]] \]
Define $L(\bv{Z}[r])$ by: 
\begin{eqnarray}
&L(\bv{Z}[r]) \defequiv \frac{1}{2}\sum_{l=1}^LZ_l[r]^2&  \label{eq:z-quad} 
\end{eqnarray} 
 Define $\Delta(\script{H}[r])$ by: 
 \[ \Delta(\script{H}[r]) \defequiv \expect{L(\bv{Z}[r+1]) - L(\bv{Z}[r])| \script{H}[r]} \]
 Suppose that per-frame changes in the  $Z_l[r]$ queues have bounded conditional fourth moments, regardless
 of past history, so that there is a constant $D>0$ such that for all $l \in \{1, \ldots, L\}$, all frames $r$, and all possible
 $\script{H}[r]$ we have: 
 \begin{equation} 
 \expect{(Z_l[r+1]- Z_l[r])^4|\script{H}[r]} \leq C \label{eq:fourth-moment} 
 \end{equation} 

\begin{thm} \label{thm:q1} (from \cite{lyap-opt2}) Suppose  $L(\bv{Z}[r])$ is the
quadratic Lyapunov function in (\ref{eq:z-quad}), and that per-frame
changes in the queues $Z_[r]$ have conditional bounded
fourth moments, so that (\ref{eq:fourth-moment}) holds. Suppose that $\expect{L(\bv{Z}[0])} < \infty$, and that 
 there are constants $B>0$, $\epsilon>0$ such that for all $r$ and all possible $\script{H}[r]$, 
 we have: 
 \begin{eqnarray*}
 &\Delta(\script{H}[r]) \leq B - \epsilon\sum_{l=1}^LZ_l[r]&
 \end{eqnarray*}
 Then for all $l \in \{1, \ldots, L\}$ we have: 
 \begin{eqnarray*}
 \sum_{r=1}^{\infty} \frac{\expect{Z_l[r]^2}}{r^2} < \infty \: \: , \: \: 
 \limsup_{R\rightarrow\infty} \frac{1}{R}\sum_{r=0}^{R-1} Z_l[r] \leq B/\epsilon \: \: (w.p.1) 
 \end{eqnarray*}
\end{thm}  

\begin{thm} \label{thm:q2} (from \cite{lyap-opt2}) Suppose  $L(\bv{Z}[r])$ is the
quadratic Lyapunov function in (\ref{eq:z-quad}), and that per-frame
changes in the queues $Z_l[r]$ have conditional bounded
fourth moments, so that (\ref{eq:fourth-moment}) holds.  Suppose that $\beta[r]$ is some additional process
related to the system, and that $\beta[0]$ and  $L(\bv{Z}[0])$ are finite with probability 1.  Suppose that:  
\[ \sum_{r=1}^{\infty} \frac{\expect{\beta[r]^2 + Z_l[r]^2}}{r^2} < \infty \]
If there are constants $V>0$, $B>0$, $\beta^*$
such that for all frames $r$ and all possible $\script{H}[r]$ 
we have: 
\begin{eqnarray}
 &\Delta(\script{H}[r]) + V\expect{\beta[r]|\script{H}[r]}  \leq B  + V\beta^*& \label{eq:dpp-cond-appc} 
 \end{eqnarray}
 Then: 
 \begin{eqnarray*} 
 \limsup_{R\rightarrow\infty} \frac{1}{R}\sum_{r=0}^{R-1}\beta[r] \leq \beta^* + B/V \: \: \: (w.p.1)
 \end{eqnarray*}
\end{thm} 

We now prove the following result. 

\begin{thm} \label{thm:appc} Suppose the same assumptions as in Theorem \ref{thm:alg1} hold, 
and for simplicity assume that $C=0$.  
Additionally assume
the fourth moment boundedness condition (\ref{eq:fourth-moment}) holds, that second moments
of $y_0[r]$ are bounded by the same constant for all $r$, and that 
there exists an $\epsilon>0$ and an  i.i.d. algorithm $\pi^*[r]$ such that: 
\begin{equation} \label{eq:slater} 
 \frac{\expect{\hat{y}_l(\pi^*[r])}}{\expect{\hat{T}(\pi^*[r])}} \leq c_l - \epsilon \: \: \forall l \in \{1, \ldots, L\} 
\end{equation} 
Then: 

(a) We have: 
 \[ \limsup_{R\rightarrow\infty} \frac{1}{R}\sum_{r=0}^{R-1} [y_0[r] - T[r]ratio^{opt}] \leq \frac{B}{V} \: \: (w.p.1) \]

(b) We have: 
\[ \limsup_{R\rightarrow\infty} \frac{\sum_{r=0}^{R-1}y_0[r]}{\sum_{r=0}^{R-1}T[r]} \leq ratio^{opt} + \frac{B}{VT^{min}} \: \: (w.p.1) \]
\end{thm} 

\begin{proof} (Theorem \ref{thm:appc} part (a)) 
Use of the history-based drift $\Delta(\script{H}[r])$ is required for the above theorems.  However, manipulations
with this drift are almost the same, and the same proof of Theorem \ref{thm:alg1} can be repeated to line
(\ref{eq:ratio-good2proof}) to show that if
the conditions of the theorem hold, then (compare with (\ref{eq:ratio-good2proof}) and note that $C=0$): 
\begin{eqnarray*}
\Delta(\script{H}[r]) + V\expect{\hat{y}_0(\pi[r])|\script{H}[r]} \leq \nonumber \\
B 
+\expect{\hat{T}(\pi[r])|\script{H}[r]}Vratio^{opt} 
\end{eqnarray*}
Rearranging terms yields: 
\begin{eqnarray}
 \Delta(\script{H}[r]) + V\expect{\hat{y}_0(\pi[r]) - \hat{T}(\pi[r])ratio^{opt}|\script{H}[r]} \leq 
B  \label{eq:dpp-plugappc} 
\end{eqnarray}
This is the same as condition (\ref{eq:dpp-cond-appc}) with $\beta[r] \defequiv \hat{y}_0(\pi[r]) - \hat{T}(\pi[r])ratio^{opt}$ and
$\beta^*\defequiv0$.   By Theorem \ref{thm:q2}, it follows that if the fourth moment 
boundedness condition (\ref{eq:fourth-moment}) holds, and if: 
\begin{equation} \label{eq:final-cond} 
 \sum_{r=1}^{\infty} \frac{\expect{\beta[r]^2 + Z_l[r]^2}}{r^2} < \infty 
 \end{equation} 
 then we can conclude the result of part (a).    
 
 It suffices to show that (\ref{eq:final-cond}) holds.  However, since we know that second moments of $y_0[r]$ and $T[r]$ are
 bounded by the same finite constant for all $r$, it is easy to see that: 
 \[ \sum_{r=1}^{\infty} \frac{\expect{\beta[r]^2}}{r^2} < \infty \]
 It suffices to show that $\sum_{r=1}^{\infty} Z_l[r]^2/r^2 < \infty$.  To this end, the proof of Theorem \ref{thm:alg1} can be repeated
 to line (\ref{eq:plug-alg1}) to show that (compare with (\ref{eq:plug-alg1})): 
 \begin{eqnarray}
&\Delta(\script{H}[r]) + V\expect{\hat{y}_0(\pi[r])|\script{H}[r]} \leq B +   \nonumber  \\
&\hspace{-.2in}\expect{\hat{T}(\pi[r])|\script{H}[r]}
 \left[ C + \frac{\expect{V\hat{y}_0(\pi^*[r]) + \sum_{l=1}^LZ_l[r]\hat{y}_l(\pi^*[r])}}{\expect{\hat{T}(\pi^*[r])}}\right] \nonumber \\
 &- \sum_{l=1}^LZ_l[r]c_l\expect{\hat{T}(\pi[r])|\script{H}[r]} \label{eq:plug-alg1b} 
\end{eqnarray}
where $\pi^*[r]$ is from any i.i.d. algorithm. 
Plugging (\ref{eq:slater})  into the right-hand-side of (\ref{eq:plug-alg1b}) yields: 
\begin{eqnarray*}
\Delta(\script{H}[r]) + V\expect{\hat{y}_0(\pi[r])|\script{H}[r]} \leq \\
B + \expect{\hat{T}(\pi[r])|\script{H}[r]}\left[C + \frac{\expect{V\hat{y}_0(\pi^*[r])}}{\expect{\hat{T}(\pi^*[r])}}\right] \\
- \expect{\hat{T}(\pi[r])|\script{H}[r]}\sum_{l=1}^LZ_l[r]\epsilon
\end{eqnarray*}
In particular, because all conditional expectations are assumed to be upper and lower bounded (and because
$T^{min}>0$), we have:  
\[ \Delta(\script{H}[r]) \leq B_1 - \epsilon T^{min}\sum_{l=1}^LZ_l[r] \]
where $B_1$ is a positive constant.
It follows by Theorem \ref{thm:q1} that $\sum_{r=1}^{\infty}\expect{Z_l[r]^2}/r^2 < \infty$ for all $l \in \{1, \ldots, L\}$.  
\end{proof}

\begin{proof} (Theorem \ref{thm:appc} part (b))
First note that if the $\limsup$ of a function $f[r]$ is upper bounded by a positive constant, then the $\limsup$ of $[f[r]]^+$ 
is bounded by the same constant, where $[f[r]]^+ \defequiv \max[f[r], 0]$.  Thus, by part (a) we have: 
\[ \limsup_{R\rightarrow\infty} \left[\frac{1}{R}\sum_{r=0}^{R-1} [y_0[r] - T[r]ratio^{opt}]\right]^+ \leq \frac{B}{V} \: \: (w.p.1) \]
Next note that because $\expect{T[r]|\script{H}[r]} \geq T^{min}>0$ for all $r$ and all $\script{H}[r]$, and second moments
of $T[r]$ are bounded by the same finite constant for all $r$, it can be shown that \cite{lyap-opt2}:
\[ \limsup_{R\rightarrow\infty} \left[\frac{R}{\sum_{r=0}^{R-1}T[r]} \right] \leq 1/T^{min} \: \: (w.p.1) \]
We then have: 
\begin{eqnarray*}
\left[\frac{\sum_{r=0}^{R-1}y_0[r]}{\sum_{r=0}^{R-1}T[r]} - ratio^{opt}\right] \leq \left[\frac{R}{\sum_{r=0}^{R-1}T[r]}\right]\times \\
\left[\frac{1}{R}\sum_{r=0}^{R-1}[y_0[r] - T[r]ratio^{opt}]  \right]^+
\end{eqnarray*}
Taking $\limsup$s of the above and using the fact that the $\limsup$ of a product of non-negative functions is less than or equal
to the product of the $\limsup$s yields: 
\[ \limsup_{R\rightarrow\infty} \left[\frac{\sum_{r=0}^{R-1}y_0[r]}{\sum_{r=0}^{R-1}T[r]} - ratio^{opt}\right] \leq \frac{1}{T^{min}}\frac{B}{V} \: \: (w.p.1) \]
\end{proof}

% ------------------------------------------------------------------------
%GATHER{Xbib.bib}   % For Gather Purpose Only
%GATHER{Paper.bbl}  % For Gather Purpose Only
\bibliographystyle{unsrt}
\bibliography{../../../latex-mit/bibliography/refs}

\begin{thebibliography}{10}

\bibitem{renewals-asilomar2010}
M.~J. Neely.
\newblock Dynamic optimization and learning for renewal systems.
\newblock {\em Proc. Asilomar Conf. on Signals, Systems, and Computers}, Nov.
  2010.

\bibitem{gallager}
R.~Gallager.
\newblock {\em Discrete Stochastic Processes}.
\newblock Kluwer Academic Publishers, Boston, 1996.

\bibitem{ross-prob}
S.~Ross.
\newblock {\em Introduction to Probability Models}.
\newblock Academic Press, 8th edition, Dec. 2002.

\bibitem{now}
L.~Georgiadis, M.~J. Neely, and L.~Tassiulas.
\newblock Resource allocation and cross-layer control in wireless networks.
\newblock {\em Foundations and Trends in Networking}, vol. 1, no. 1, pp. 1-149,
  2006.

\bibitem{neely-mdp-cdc09}
M.~J. Neely.
\newblock Stochastic optimization for markov modulated networks with
  application to delay constrained wireless scheduling.
\newblock {\em Proc. IEEE Conf. on Decision and Control (CDC)}, Shanghai,
  China, Dec. 2009.

\bibitem{self-learning-mdp}
F.~J.~V\'{a}zquez Abad and V.~Krishnamurthy.
\newblock Policy gradient stochastic approximation algorithms for adaptive
  control of constrained time varying markov decision processes.
\newblock {\em Proc. IEEE Conf. on Decision and Control}, Dec. 2003.

\bibitem{q-learning-mimo}
D.~V. Djonin and V.~Krishnamurthy.
\newblock $q$-learning algorithms for constrained markov decision processes
  with randomized monotone policies: Application to mimo transmission control.
\newblock {\em IEEE Transactions on Signal Processing}, vol. 55, no. 5, pp.
  2170-2181, May 2007.

\bibitem{energy-delay-approx}
N.~Salodkar, A.~Bhorkar, A.~Karandikar, and V.~S. Borkar.
\newblock An on-line learning algorithm for energy efficient delay constrained
  scheduling over a fading channel.
\newblock {\em IEEE Journal on Selected Areas in Communications}, vol. 26, no.
  4, pp. 732-742, May 2008.

\bibitem{mihaela-dp1}
F.~Fu and M.~van~der Schaar.
\newblock A systematic framework for dynamically optimizing multi-user video
  transmission.
\newblock {\em IEEE Journal on Selected Areas in Communications}, vol. 28, no.
  3, pp. 308-320, April 2010.

\bibitem{mihaela-dp2}
F.~Fu and M.~van~der Schaar.
\newblock Decomposition principles and online learning in cross-layer
  optimization for delay-sensitive applications.
\newblock {\em IEEE Trans. Signal Processing}, vol. 58, no. 3, pp. 1401-1415,
  March 2010.

\bibitem{neely-thesis}
M.~J. Neely.
\newblock {\em Dynamic Power Allocation and Routing for Satellite and Wireless
  Networks with Time Varying Channels}.
\newblock PhD thesis, Massachusetts Institute of Technology, LIDS, 2003.

\bibitem{neely-energy-it}
M.~J. Neely.
\newblock Energy optimal control for time varying wireless networks.
\newblock {\em IEEE Transactions on Information Theory}, vol. 52, no. 7, pp.
  2915-2934, July 2006.

\bibitem{neely-fairness-infocom05}
M.~J. Neely, E.~Modiano, and C.~Li.
\newblock Fairness and optimal stochastic control for heterogeneous networks.
\newblock {\em Proc. IEEE INFOCOM}, March 2005.

\bibitem{sno-text}
M.~J. Neely.
\newblock {\em Stochastic Network Optimization with Application to
  Communication and Queueing Systems}.
\newblock Morgan \& Claypool, 2010.

\bibitem{atilla-fairness-ton}
A.~Eryilmaz and R.~Srikant.
\newblock Fair resource allocation in wireless networks using
  queue-length-based scheduling and congestion control.
\newblock {\em IEEE/ACM Transactions on Networking}, vol. 15, no. 6, pp.
  1333-1344, Dec. 2007.

\bibitem{stolyar-greedy}
A.~Stolyar.
\newblock Maximizing queueing network utility subject to stability: Greedy
  primal-dual algorithm.
\newblock {\em Queueing Systems}, vol. 50, no. 4, pp. 401-457, 2005.

\bibitem{stolyar-gpd-gen}
A.~Stolyar.
\newblock Greedy primal-dual algorithm for dynamic resource allocation in
  complex networks.
\newblock {\em Queueing Systems}, vol. 54, no. 3, pp. 203-220, 2006.

\bibitem{primal-dual-cmu}
Q.~Li and R.~Negi.
\newblock Scheduling in wireless networks under uncertainties: A greedy
  primal-dual approach.
\newblock {\em Arxiv Technical Report: arXiv:1001:2050v2}, June 2010.

\bibitem{lin-shroff-cdc04}
X.~Lin and N.~B. Shroff.
\newblock Joint rate control and scheduling in multihop wireless networks.
\newblock {\em Proc. of 43rd IEEE Conf. on Decision and Control, Paradise
  Island, Bahamas}, Dec. 2004.

\bibitem{vijay-allerton02}
R.~Agrawal and V.~Subramanian.
\newblock Optimality of certain channel aware scheduling policies.
\newblock {\em Proc. 40th Annual Allerton Conference on Communication ,
  Control, and Computing, Monticello, IL}, Oct. 2002.

\bibitem{prop-fair-down}
H.~Kushner and P.~Whiting.
\newblock Asymptotic properties of proportional-fair sharing algorithms.
\newblock {\em Proc. of 40th Annual Allerton Conf. on Communication, Control,
  and Computing}, 2002.

\bibitem{chihping-utility-round-robin}
C.-P. Li and M.~J. Neely.
\newblock Network utility maximization over partially observable markovian
  channels.
\newblock {\em Arxiv Technical Report: arXiv:1008.3421v1}, Aug. 2010.

\bibitem{network-corroboration-arxiv}
B.~Liu, P.~Terlecky, A.~Bar-Noy, R.~Govindan, and M.~J. Neely.
\newblock Optimizing information credibility in social swarming applications.
\newblock {\em ArXiv technical report, arXiv:1009:6006}, Sept. 2010.

\bibitem{williams-martingale}
D.~Williams.
\newblock {\em Probability with Martingales}.
\newblock Cambridge Mathematical Textbooks, Cambridge University Press, 1991.

\bibitem{lyap-opt2}
M.~J. Neely.
\newblock Queue stability and probability 1 convergence via lyapunov
  optimization.
\newblock {\em Arxiv Technical Report}, Oct. 2010.

\bibitem{bertsekas-neural}
D.~P. Bertsekas and J.~N. Tsitsiklis.
\newblock {\em Neuro-Dynamic Programming}.
\newblock Athena Scientific, Belmont, Mass, 1996.

\bibitem{neely-mwl-arxiv}
M.~J. Neely.
\newblock Max weight learning algorithms with application to scheduling in
  unknown environments.
\newblock {\em arXiv:0902.0630v1}, Feb. 2009.

\bibitem{bisdikian-qoi}
C.~Bisdikian, L.~M. Kaplan, M.~B. Srivastava, D.~J. Thornley, D.~Verma, and
  R.~I. Young.
\newblock Building principles for a quality of information specification for
  sensor information.
\newblock {\em 12th Int'l Conf. on Information Fusion (Fusion '09), Seattle,
  WA}, July 2009.

\end{thebibliography}
\end{document}